\documentclass[11pt]{article}
\usepackage[latin1]{inputenc}
\usepackage{amsmath,amssymb}
\usepackage{color}
\usepackage{multicol}
\usepackage{floatrow}
\usepackage{latexsym,epsfig}
\usepackage{amssymb,amsmath,amsfonts,amscd}
\usepackage{exscale}
\usepackage{ifthen}
\usepackage{longtable}
\usepackage{subfigure}
\usepackage{cite}
\usepackage{lscape}
\usepackage{blindtext}

\newtheorem{theo}{\textbf{Theorem}\ }
[section]
\newtheorem{lem}[theo]{\textbf{Lemma}\ }

\newtheorem{prop}[theo]{\textbf{Proposition}\ }

\def \N{\mathbb{N}}

\def \E{\mathbb{E}}
\def \X{\mathbb{X}}
\def \R{\mathbb{R}}

\def \P{\mathcal{P}}

	\definecolor{amber}{rgb}{1.0, 0.75, 0.0}
	\definecolor{airforceblue}{rgb}{0.36, 0.54, 0.66}
		\definecolor{alizarin}{rgb}{0.82, 0.1, 0.26}
	\definecolor{amaranth}{rgb}{0.9, 0.17, 0.31}	
	\definecolor{amber(sae/ece)}{rgb}{1.0, 0.49, 0.0}	
	\definecolor{americanrose}{rgb}{1.0, 0.01, 0.24}	
	\definecolor{amethyst}{rgb}{0.6, 0.4, 0.8}	
	\definecolor{antiquebrass}{rgb}{0.8, 0.58, 0.46}	
	\definecolor{ao(english)}{rgb}{0.0, 0.5, 0.0}	
	\definecolor{apricot}{rgb}{0.98, 0.81, 0.69}	
	\definecolor{aquamarine}{rgb}{0.5, 1.0, 0.83}	
	\definecolor{armygreen}{rgb}{0.29, 0.33, 0.13}
	\definecolor{arsenic}{rgb}{0.23, 0.27, 0.29}	
	\definecolor{asparagus}{rgb}{0.53, 0.66, 0.42}	
	\definecolor{auburn}{rgb}{0.43, 0.21, 0.1}	
	\definecolor{awesome}{rgb}{1.0, 0.13, 0.32}	
	\definecolor{ballblue}{rgb}{0.13, 0.67, 0.8}	
	\definecolor{bananayellow}{rgb}{1.0, 0.88, 0.21}	
	\definecolor{bittersweet}{rgb}{1.0, 0.44, 0.37}	
	\definecolor{bleudefrance}{rgb}{0.19, 0.55, 0.91}
	\definecolor{blush}{rgb}{0.87, 0.36, 0.51}	
	\definecolor{bostonuniversityred}{rgb}{0.8, 0.0, 0.0}	
	\definecolor{brickred}{rgb}{0.8, 0.25, 0.33}	
\definecolor{brightmaroon}{rgb}{0.76, 0.13, 0.28}	
	\definecolor{brilliantrose}{rgb}{1.0, 0.33, 0.64}
	\definecolor{britishracinggreen}{rgb}{0.0, 0.26, 0.15}	
	\definecolor{bulgarianrose}{rgb}{0.28, 0.02, 0.03}	
	\definecolor{burgundy}{rgb}{0.5, 0.0, 0.13}	
	\definecolor{burntorange}{rgb}{0.8, 0.33, 0.0}	
	\definecolor{byzantium}{rgb}{0.44, 0.16, 0.39}			
	\definecolor{candyapplered}{rgb}{1.0, 0.03, 0.0}
	\definecolor{capri}{rgb}{0.0, 0.75, 1.0}	
	\definecolor{carminered}{rgb}{1.0, 0.0, 0.22}	
	\definecolor{charcoal}{rgb}{0.21, 0.27, 0.31}	
	\definecolor{coquelicot}{rgb}{1.0, 0.22, 0.0}	
	\definecolor{daffodil}{rgb}{1.0, 1.0, 0.19}

\numberwithin{equation}{section}

\begin{document}

\begin{center}
 {\huge \bf Conditioned limit theorems for products of positive random matrices } 

\vspace{1cm}

  C. Pham 
$^($\footnote{
Universit\'e Fr. Rabelais Tours,  LMPT  UMR  CNRS 7350,  Tours,  France. \\ 	email : Thi-Da-Cam.Pham@lmpt.univ-tours.fr}$^)$
$^($\footnote{ The author thanks the {\bf Vietnam Institute for Advanced Studies in Mathematics} (VIASM) in Ha Noi for generous hospitality and financial support during the first semester 2017. }$^)$

\vspace{0.5cm}
     \end{center}
     
        \centerline{\bf  \small Abstract }

        Inspired by  a recent paper of  I. Grama, E. Le Page and M. Peign\'e, we consider a sequence $(g_n)_{n \geq 1}$ of i.i.d. random $d\times d$-matrices with non-negative entries and study the fluctuations of the process $(\log \vert g_n\cdots g_1\cdot x\vert )_{n \geq 1}$  for any  non-zero vector $x$ in $\mathbb R^d$ with non-negative coordinates. Our method involves approximating  this process by a martingale and  studying   harmonic functions for its restriction  to the upper half line.  Under certain conditions, the probability for this process to stay in the upper half real line up to time $n$ decreases as $c \over \sqrt n$ for some positive constant $c$.

       \vspace{0.5cm}



\section{Introduction}
 
\hspace{0.5cm }   Many limit theorems  describe the asymptotic behaviour of random walks with i.i.d. increments, for instance the strong law of large numbers,  the central limit theorem,  the invariant principle... Besides, the   fluctuations of these processes  are well studied, for example the decay of the probability  that they stay inside the half real line up to time $n$.   The Wiener-Hopf factorization is usually used  in this case but so far,   it seems to be impossible to adapt in   non-abelian context. Recently,  much efforts are made  to apply the results above for the  logarithm of the norm of the product of i.i.d. random matrices since  it behaves  similarly to  a sum of i.i.d. random variables. Many limit theorems arose for the last 60 years, initiated by Furstenberg-Kersten \cite{FK}, Guivarc'h-Raugi \cite{GR}, Le Page \cite{LePage82}... and recently Benoist-Quint \cite{BQ}. Let us mention also the works by  Hennion \cite{H1} and Hennion-Herve \cite{HH} for matrices with positive entries. However, the studies on the subject of fluctuation was quite sparse a few years ago. Thanks to the approach of Denisov and Wachtel  \cite{DW} for random walks in  Euclidean spaces and motivated by branching processes, I. Grama, E. Le Page and M. Peign\'e  recently progressed in \cite{GLP1} for invertible matrices.  Here we propose to develop the same strategy for matrices with positive entries by using \cite{HH}.


\noindent  We endow $\mathbb R^d$ with the norm $\vert \cdot \vert $ defined by   $\displaystyle  \vert x\vert := \sum_{i=1}^d \vert x_i\vert $ for any column vector $x=(x_i)_{1\leq i \leq d}$.
 Let $\mathcal C$ be the cone of vectors in $\mathbb R^d$ with non-negative coordinates 
 $$
 \mathcal C  := \{x \in \mathbb R ^d: \forall 1 \le i \le d, x_i \ge 0 \}
 $$ 
 and   $\mathbb X$ be the limited cone defined by $$
 \mathbb X := \{ x \in \mathcal C,  |x| =1 \}.
 $$
 
  \noindent Let $S $ be the set of $d\times d$ matrices with non-negative entries such that  each column contains at least one  positive entry; its interior is $\mathring {S }:=\{g =(g(i, j))_{1\leq i, j \leq d}/ g(i, j) >0\} $. Endowed with the standard multiplication of matrices, the set $S$ is a semigroup  and $ \mathring S$ is the ideal of $S$, more precisely, for any $g \in \mathring S$ and  $h \in  S$, it is evident that $gh \in \mathring S$.

 \noindent We consider the following actions:
 
 \begin{itemize}
 	\item the left  linear action of $S$ on $\mathcal C$ defined by $(g,x) \mapsto gx$ 	for any  $g \in S$ and $x \in \mathcal C$, 
 	
 	\item the left projective action of $S$ on $\mathbb X$ defined by $(g,x) \mapsto g \cdot x := \frac{gx}{|gx|}$ for any  $g \in S$ and $x \in \mathbb X$.
 \end{itemize}
 For any $g =(g(i, j))_{1 \le i, j \le d} \in S$,  without confusion, let   
 \vspace{-0.1cm} 
 $$v(g) := \min_{1\leq j\leq d}\Bigl(\sum_{i=1}^d g(i, j)\Bigr)\quad {\rm and} \quad  \vert g\vert :=\max_{1\leq j\leq d}\Bigl(\sum_{i=1}^d g(i, j)\Bigr) , 
 $$ 
 then $\vert \cdot \vert$ is a norm on $S$ and for any $x \in \mathcal C$,
 \begin{equation}\label{controlnormgx}
 0< v(g)\  \vert x\vert \leq \vert gx\vert \leq \vert g\vert \ \vert x\vert.
 \end{equation}
 We set $N(g):= \max \left({1\over v(g)}, \vert g\vert\right) $.

 On some probability space $(\Omega, \mathcal F, \mathbb P)$, we consider a sequence of i.i.d.  $S$-valued matrices $ (g_n)_{n \geq 0}$ with the same distribution  $\mu$   on $S$. Let  $L_0= Id$ and $  L_n := g_n \ldots g_1$ for any $n \ge 0$. For any fixed $x \in \mathbb X$, we define  the $\mathbb X$-valued Markov chain $(X_n^x)_{n \geq0}$ by setting $X_n^x := L_n \cdot x$ for any $n \ge 0$ (or simply $X_n$ if there is no confusion). We denote by $P$ the transition probability of $(X_n)_{n \ge 0}$,  defined by: for any $  x  \in \mathbb{X}$ and any bounded Borel function $\varphi :  \mathbb X \to \mathbb{C}$, 
 \vspace{-0.1cm}
 \begin{eqnarray*} 
 P\varphi ( x) := \int_{S} {\varphi ( g\cdot x)  \mu (dg)} =\mathbb E[\varphi(L_1\cdot x)].
 \end{eqnarray*}
 Hence, for any  $n \ge 1$,
  $$ P^n\varphi(x) =\mathbb E[\varphi(L_n\cdot x)].$$

 We assume that with positive probability, after finitely many  steps, the sequence $(L_n)_{n\ge 1}$ reaches $\mathring S$. In mathematical term, it is equivalent to writing as
 \vspace{-0.1cm}
\begin{eqnarray*}
\mathbb P \left( \bigcup _{n \ge 1} [L_n \in \mathring S] \right) >0.
\end{eqnarray*}
On the product space $S \times \mathbb X$, we define the function $\rho$ by setting for any $(g,x) \in S \times \mathbb X $,
 $$
 \rho(g,x):=\log |gx|.
 $$
 Notice that $gx = e^{\rho(g,x)} g \cdot x$; in other terms, the linear action of $S$ on $\mathcal C$ corresponds to the couple $(g \cdot x, \rho(g,x))$. This function $\rho$ satisfies the cocycle property  $\rho(gh, x)=  \rho(g, h\cdot x)+\rho(h, x)$ for any $g, h\in S$ and $x \in \mathbb  X$ and implies the basic decomposition for any $x \in \mathbb X$,
 \begin{equation*} 
 \log \vert L_n x \vert =  \sum_{k=1}^n  \rho(g_k, X^x_{k-1}).
 \end{equation*} 
For any $a \in \mathbb R$ and $n \ge 1$, let $S_0 := a$ and $S_n = S_n (x,a) := a + \sum_{k=1}^n  \rho(g_k, X_{k-1})$.
Then the sequence  $(X_n, S_n)_{n \geq0}$ is a Markov chain on $\mathbb X \times \mathbb R$ with transition  probability $\widetilde P$ defined by: for any $(x, a) \in \mathbb X \times \mathbb R$ 
  and any bounded Borel function $\psi : \mathbb X \times \mathbb R \to \mathbb{C}$, 
 $$\widetilde P \psi ( x, a) = \int_{S} {\psi (g \cdot x, a+\rho( g, x)) \mu (dg)} .$$
  For any $(x, a)\in \mathbb{X}\times \mathbb R$, we denote by   $\mathbb{P}_{ x, a}$  the probability measure  on $(\Omega, \mathcal F)$  conditioned to the event $
 [X_{0}=x, S_0 =  a]$  and by $\mathbb{E}_{x, a}$ 
 the corresponding expectation; for the sake of brevity, by $\mathbb{P}_{ x}$  we denote $\mathbb P _{x,a}$ when $S_0 =0$  and by $\mathbb{E}_{x}$  the corresponding expectation. 
 Hence for any  $n \ge 1$,
 \begin{eqnarray*} \label{eqn4}
 \widetilde P^n\psi(x, a) = \mathbb E[\psi(L_n\cdot x, a+\log \vert L_n x\vert)]= \mathbb E_{x,a} [\psi( X_n,  S_n)].
 \end{eqnarray*}
Now we consider the restriction of $\widetilde P_+$ to $\mathbb X \times \mathbb R^+$ of $\widetilde P$ defined by: for any $(x,a) \in \mathbb X \times \mathbb R ^+$ and any bounded function $\psi : \mathbb X \times \mathbb R \to \mathbb C$, 
 $$
 \widetilde P _+ \psi (x, a) =  P (\psi {\bf 1}_{\mathbb X \times \mathbb R^+_\ast})(x,a).
 $$
 Let us emphasize that $\widetilde P_+$ may not be a Markov kernel on $\mathbb X \times \mathbb R ^+$.

  Let $  \tau := \min \{ n \ge 1:\, {S_n} \leq 0\}  $ be  the first time the random process $(S_n)_{n \ge 1}$ becomes non-positive; for any $(x, a) \in \mathbb X \times \mathbb{R^+}$ and any bounded Borel function $\psi : \mathbb X \times \mathbb R \to \mathbb{C}$, 
\begin{eqnarray*}\label{eqn5}
\widetilde P_+\psi(x, a)=\mathbb E_{x, a}[\psi(X_1, S_1); \tau >1]= \mathbb E [\psi (g_1\cdot x, a+ \rho(g_1,x));   a+\rho(g_1,x) > 0].
\end{eqnarray*}

A positive $\widetilde P_+$-harmonic function $V$ is any function from $\mathbb X \times \mathbb R ^+$ to $\mathbb R ^+$  satisfying $ \widetilde P_+ V =V $. We extend $V$ by setting $V(x, a) =0$ for $(x, a) \in \mathbb X \times \mathbb R^-_\ast$.

 In other words, the function $V$ is $\widetilde P_+$-harmonic if and only if for any $x \in \mathbb X$ and  $a\ge0$,
 \begin{eqnarray} \label{eqn7}
  V(x, a) = \mathbb E_{x,a} [V(X_1, S_1); \tau >1].
 \end{eqnarray}

 
 From Theorem II.1 in \cite{HH}, under conditions P1-P3 introduced below, there exists a unique probability measure $\nu$ on $\mathbb X$ such that for any bounded Borel function $\varphi$ from $\mathbb X$ to $\mathbb R$,
 $$
 (\mu \ast \nu) (\varphi)= \int_{S} \int_{\mathbb X} \varphi(g \cdot x) \nu(dx) \mu(dg) = \int_{\mathbb X} \varphi(x) \nu(dx) = \nu (\varphi).
 $$
 Such a measure is said to be $\mu$-invariant.
 Moreover, the upper Lyapunov exponent associated with $\mu$ is finite and is expressed by
 \begin{eqnarray} \label{eqn10}
 \gamma_\mu = \int_{S} \int_{\mathbb X}  \rho(g,x) \nu (dx) \mu (dg) .
 \end{eqnarray}
 
 Now we state the needed hypotheses for later work.
 
  \noindent {\bf HYPOTHESES}

 \noindent {\bf P1} {\it  There exists $\delta_0>0$ such that $\displaystyle \int_{S}   N(g)^{ \delta_0} \mu(dg) <+\infty$.}

 \noindent {\bf P2}  {\it There exists no affine subspaces $A$ of $\R ^d$ such that $A \cap \mathcal C$ is non-empty and bounded and invariant under the action of all elements of the support of $\mu$.}

 \noindent {\bf P3}   {\it There exists $n_0 \geq 1$ such that  $ \mu^{*n_0}(\mathring {S }) >0.$}
 

 \noindent {\bf P4} {\it The upper Lyapunov exponent $\gamma_\mu$ is equal to $0$.}
 
 \noindent {\bf P5} {\it There exists  $\delta >0$  such that $\mu \{g \in S : \forall x \in \mathbb X, \log|gx| \ge \delta \} >0$.}


In this paper, we establish the asymptotic behaviour of $\mathbb P_{x,a} (\tau >n)$ by studying the $\widetilde P_+$-harmonic function $V$. More precisely,  Proposition \ref{theo2}  concerns the existence of a $P^+$-harmonic function and its properties whereas Theorem \ref{theo3} is about the limit behaviour of the exit time $\tau$.
 \begin{prop} \label{theo2}
 	Assume hypotheses P1-P5.
 	\begin{enumerate}
 		\item  For any $x \in \mathbb X$ and $a \ge 0$, the sequence $\Bigl(\mathbb E _{x,a} [ S_n; \tau  >n ] \Bigr) _{n \ge 0}$ converges to the function $	V(x, a):= a- \mathbb E_{x,a} M_\tau$.
 		
 			\item For  any $x \in \mathbb X$ the function $V(x, \cdot)$ is increasing on $\mathbb R^+ $.
 		
 		\item There exists $c >0$ and $A >0$ such that for any $x \in \mathbb  X$ and $a\ge 0$,  
 		$$
 		\frac{1}{c} \vee (a-A) \leq V(x, a) \leq c(1+a).
 		$$
 	 		
 		\item For any $x \in \mathbb X$, the function $V(x, .)$ satisfies $\displaystyle \lim_{a \to +\infty} \frac{V(x, a)}{a} =1$.

 		\item The function $V$ is $\widetilde P_+$-harmonic.
 		
 	\end{enumerate}
 \end{prop}
 
 The function $V$ contains information of the part of the trajectory which stays in $\mathbb R^+$ as stated in Theorem \ref{theo3}.
 
 \begin{theo} \label{theo3}
 	Assume P1-P5. Then for any $x \in \mathbb X$ and $a \ge0$,
 	$$
 	\mathbb P _{x,a} (\tau >n) \sim \frac{2V(x, a)}{\sigma \sqrt{2 \pi n}} \,\,\mbox{as} \,\, n \to +\infty.
 	$$
 	Moreover, there exists a constant $c$ such that for any $x \in \mathbb X$, $a\ge 0$ and $n \ge 1$,
 	$$
  \sqrt n \mathbb P _{x,a} (\tau  > n) \leq c V(x,a).
 	$$

 \end{theo}

 As a direct consequence, we prove that the sequence $(\frac{S_n}{\sigma \sqrt n})_{n \ge 1}$, conditioned to the event $\tau >n$, converges in distribution to the Rayleigh law as stated below. 
 
 \begin{theo} \label{theo4} 
 	Assume P1-P5. For any $x \in \mathbb X$, $a\ge 0$ and $t >0$, 
 	$$
 	\lim_{n \to +\infty} \mathbb P _{x,a} \left( \frac{S_n}{ \sqrt n} \leq t \mid \tau >n\right)  = 1- \exp \left(- \frac{t^2}{2\sigma ^2} \right).
 	$$
 
 \end{theo}

 In section 2, we approximate the chain $(S_n)_{n \ge 0}$ by a martigale and in  section 3, we study the harmonic function $V$ and state the proof of Proposition \ref{theo2}. We use the coupling argument to prove Theorem \ref{theo3} and Theorem \ref{theo4} in section 4. At last, in section 5 we check general conditions to apply an invariant principle stated in Theorem 2.1 in \cite{GLP1}.

 Throughout this paper, we denote the absolute constants such as $C, c,c_1, c_2, \ldots$ and the constants depending on their indices such as $c_\varepsilon, c_p, \ldots$. Notice that they are not always the same when used in different formulas. The integer part of a real constant $a$ is denoted by $[a]$.
 

  \section{Approximation of the chain $(S_n)_{n \ge 0}$} 
   
   In this section, we discuss  the spectral properties of $P$ and then utilise them to approximate the chain $(S_n)_{n \ge 0}$. Throughout this section, we assume that conditions P1-P4 hold true.
  \subsection{Spectral properties of the operators $P$ and its Fourier transform}

 Following \cite{H2}, we endow $\mathbb X$ with a bounded distance $d$ such that $g$ acts on $\mathbb X$ as a contraction with respect to $d$ for any $g \in S$.  For any $x,y \in \mathbb X$, we write: 
 \[\begin{array}{l}
 \displaystyle m\left( {x,y} \right) =  \min_{1 \le i \le d} \left\{ {\displaystyle \frac{{{x_i}}}{{{y_i}}} |  {y_{i } >0}} \right\}
 \end{array}\]
 and it is clear that $ 0 \le m\left( {x,y} \right) \le 1$. For any $x,y \in \mathbb X$, let $ d\left( {x,y} \right): = \varphi \left( {m\left( {x,y} \right)m\left( {y,x} \right)} \right)$, where  $\varphi$ is the one-to-one function  defined for any $s \in [0,1]$ by $\varphi \left( s \right): = \displaystyle \frac{1 - s}{1 + s}$. Set \(c\left( g \right): = \sup \left\{ {d\left( {g\cdot x,g\cdot y} \right),x,y \in \mathbb{X}} \right\}\) for $g \in S $; the proposition below gives some more properties of $d$ and $c(g)$.
 
 \begin{prop} \cite{H2} \label{propH}
 	 The quantity $d$ is a distance on ${\mathbb{X}}$ satisfying the following properties: 
 	\begin{enumerate}
 		\item $\sup \{  d(x,y ): x, y \in  {\mathbb{X}} \} =1 $.
 		\item $|x-y| \le 2d(x,y)$ for any $x,y \in \mathbb X$. 
 		\item $c(g) \le 1$ for any $g \in S$, and $c(g) < 1$ if and only if $g \in \mathring S$.
 		\item \(d\left( {g\cdot x,g\cdot y} \right) \le c\left( g \right)d\left( {x,y} \right) \le c(g)\) for any  and \(x,y \in {\mathbb{X}}\).
 		\item \(c\left( {gh} \right) \le c\left( g \right)c\left( {h} \right)\) for any \(g, h \in S \).
 	\end{enumerate}	
 	
 \end{prop}
 
 From now on, we consider a sequence $(g_n)_{n \ge 0}$ of i.i.d. $S$-valued random variables, we set $a_k := \rho (g_k, X_{k-1})$ for $k \ge 1$ and hence $S_n = a + \sum_{k=1}^n a_k$ for $n \ge 1$. In order to study the asymptotic behavior of the process $(S_n)_{n \geq 0}$, we need to consider the ``Fourier transform''  of the random variables $a_k$, under $\mathbb P_x, x \in \mathbb X$, similarly for classical random walks with independent increments  on $\mathbb R$.  Let  $P_t$ be the family of ``Fourier operators'' defined for any $t \in \mathbb R$,  $  x  \in \mathbb{X}$ and any bounded Borel function $\varphi :  \mathbb X \to \mathbb{C}$ by
 \begin{eqnarray} \label{eqn201}
 P_t\varphi ( x) := \int_{S} e^{it \rho(g, x)}\varphi ( g\cdot x)  \mu (dg)= \mathbb E_x\left[ e^{it a_1} \varphi(X_1) \right]
 \end{eqnarray}
 and for any $ n \geq 1$,
 \begin{eqnarray} \label{eqn8}
 P^n_t\varphi ( x)= \mathbb E[e^{it \log \vert L_n x\vert} \varphi(L_n\cdot x)]=\mathbb E_x [e^{it S_n} \varphi(X_n)].
 \end{eqnarray}
 Moreover,  we can imply that
 \begin{eqnarray}\label{eqn1.131}
 P^mP_t^n\varphi(x)&=&\mathbb E \left[ e^{it\log\vert g_{m+n}\cdots g_{m+1}(L_m\cdot x)\vert}\varphi(L_{m+n}\cdot x)\right] \notag \\
 &=& \mathbb E_x \left[ e^{it (a_{m+1}+\cdots +a_{m+n})}\varphi(X_{n+m}) \right]
 \end{eqnarray}
 and when $\varphi = 1$, we obtain
 \begin{eqnarray*} \label{eqn91}
 \mathbb E_x \left[ e^{it S_n}\right] = P^n_t1(x)\quad {\rm and} \quad \mathbb E_x \left[ e^{it (a_{m+1}+\cdots +a_{m+n})} \right] = P^mP_t^n 1(x).
 \end{eqnarray*}

 We consider the space $C(\mathbb X)$ of continuous functions from $\mathbb X$ to $\mathbb C$ endowed with the norm of uniform convergence $\vert.\vert_\infty$. Let $L$ be the subset of Lipschitz functions on $\mathbb X$ defined by 
 $$
 L := \{\varphi \in C(\mathbb X): |\varphi|_L := |\varphi|_\infty + m(\varphi)   <+\infty\},
 $$
where $m(\varphi):= \sup_{\stackrel{x, y \in \mathbb X}{x\neq y}} {\vert \varphi(x)-\varphi(y)\vert \over d(x, y)}$.
 The spaces $(C(\mathbb X), \vert \cdot\vert_\infty)$ and $(L, \vert \cdot\vert _L)$ are  Banach spaces and  the canonical injection from $L$ into $C(\mathbb X)$ is compact. The norm of a bounded operation $A$ from $L$ to $L$ is denoted by $|A|_{L \to L} := \sup_{\varphi \in L} |A \varphi|_L$. We denote $L'$ the topological dual of $L$ endowed  with the norm $\vert\cdot \vert_{L'}$ corresponding to $\vert \cdot \vert _L$; notice that any probability measure $\nu$ on $\mathbb X$ belongs to $L'$. 
 
 \noindent For further uses, we state here some helpful estimations. 

\begin{lem} \label{lemHH}
	For $g \in S$, $x,y,z \in \mathbb X$ such that $d(x,y) <1$ and for any $t\in \mathbb R$, 
	\begin{eqnarray} \label{lem2.21}
	|\rho(g,x)| \le 2 \log N(g),
	\end{eqnarray}
	and 
	\begin{eqnarray} \label{lem2.22}
	|e^{it\rho(g,y)} - e^{it\rho(g,z)} | \le \Bigl( 4 \min\{2|t| \log N(g), 1 \}  + 2 C |t| \Bigr) d(y,z),
	\end{eqnarray}
	where $C= \sup \{ \frac{1}{u} \log \frac{1}{1-u} : 0 < u \le \frac{1}{2} \} <+\infty$.

\end{lem}
{\bf Proof.}
For the first assertion, from (\ref{controlnormgx}), we can imply that $\vert \log |gx| \vert \le \log N(g)$.
For the second assertion, we refer to the proof the Theorem III.2 in \cite{HH}.

\rightline{$\square$}

Denote  $\varepsilon(t) := \int_S \min \{ 2 |t| \log N(g) ,2\} \mu (dg)$. Notice that $\lim_{t \to 0} \varepsilon(t) =0$.

 \begin{prop} \cite{HH} \label{spectre} Under hypotheses {\bf P1, P2, P3} and {\bf P4}, for any $t \in \mathbb R$, the operator $P_t$ acts on $L$ and satisfies the following properties:
 	\begin{enumerate}
 		\item Let $\Pi : L \to L$ be the rank one operator   defined by $\Pi(\varphi)= \nu(\varphi ) 1$ for any function $\varphi \in L$,  where $\nu$ is the unique $P$-invariant probability measure on $\mathbb X$ and  $R:= P- \Pi$.
 		
 		The operator $R: L\to L$ satisfies 
 			\begin{equation*} 
 		 \Pi R=R\Pi=0 ,
 		\end{equation*}
 		and its spectral radius is less than $1$; in other words,  there exist  constants $ C>0$  and $0 < \kappa < 1$ such that  $\vert R ^n\vert_{L \to L}\leq C \kappa^n$ for any $n \geq 1$.
 		 		
 		\item There exist  $\epsilon>0$ and $ 0 \le r_\epsilon <1$ such that for any $t\in [-\epsilon, \epsilon]$, there exist a complex number $\lambda_t$ closed to $1$   with modulus less than or equal to $ 1$, a rank one operator $\Pi_t$  and an operator $R_t$  on $L$   with spectral radius less than or equal to $ r_\epsilon$ such that
 		$$
 		P_t= \lambda_t \Pi_t+R_t  \quad {\it and} \quad \Pi_t R_t=R_t \Pi_t = 0.
 		$$
 		
 	Moreover,  $C_P:= \displaystyle  \sup_{\stackrel{-\epsilon\le t\leq \epsilon}{n \geq 0}}\vert P_t^n\vert _{L \to L}<+\infty$.
 	
 		\item For any $p\geq 1$,
 		\begin{equation}\label{momentsoforderp}
 		\sup_{n \geq  0}\sup_{x \in \mathbb X}\mathbb E_x \vert \rho(g_{n+1}, X_n)\vert ^p<+\infty.
 		\end{equation}
 		
 	\end{enumerate}
 \end{prop}

 \noindent{ \bf Proof.}
 {\bf (a)} We first check that  $P_t$ acts on $(L, \vert\cdot \vert _L)$ for any $t \in \mathbb R$. On one hand,  $|P _t \varphi|_\infty \leq |\varphi|_\infty $ for any  $\varphi \in L$. On the other hand, by (\ref{lem2.22}) for any $x, y \in \mathbb X$ such that $ x \ne y$, 
 \begin{eqnarray*}
 	\frac {|P _t \varphi (x)-  P _t \varphi (y)| }{d(x,y)}
 	&\le& \int_{S} \left( \left| \frac {e^{it\rho(g,x)} - e^{it\rho(g,y)} } {d(x,y)} \right| |\varphi(g \cdot x)| +  \left| \frac {\varphi(g \cdot x) -\varphi(g \cdot y)} {d(x,y)} \right|\right) \mu(dg) \\
 	&\le& |\varphi|_\infty (4 \varepsilon(t) + 2C |t| )    +  \int_{S} \left(  \frac{| \varphi(g \cdot x) -\varphi(g \cdot y)|}{d(g \cdot x, g \cdot y)}  \frac{d(g \cdot x, g \cdot y)}{d(x,y)} \right) \mu(dg) ,\\
 	&\le&  |\varphi|_\infty (4 \varepsilon(t) + 2C |t| )    + m(\varphi),
 \end{eqnarray*}
 
 \noindent which implies $m(P _t \varphi)   \leq |\varphi|_\infty (4 \varepsilon(t) + 2C |t| )    + m(\varphi) < +\infty$. Therefore $P _t \varphi \in L$.
 
 \vspace{0.3cm}
 {\bf (b)} Let $\Pi$ be the rank one projection on $L$ defined by $\Pi \varphi = \nu(\varphi) {\bf 1}$ for any $\varphi \in L$. Let $ R := P - \Pi$. By definition, we obtain $P \Pi = \Pi P = \Pi$ and $\Pi ^2 = \Pi$  which implies $\Pi R = R \Pi =0$ and $R^n = P^n - \Pi$ for any $n \ge 1$. Here we only sketch the main steps by taking into account the proof of  Theorem III.1 in \cite{HH}. 
 
  \noindent Let $\mu ^{*n}$ be the distribution of the random variable $L_n$ and set
 \begin{eqnarray*} \label{eqn16}
 c(\mu^{*n}) := \sup \left\{  \int _S \frac{d(g \cdot x, g \cdot y)}{d(x,y)} d \mu^{*n}(g) : x,y \in \mathbb X, x \ne y \right\}.
 \end{eqnarray*}
 Since $c(\cdot) \le 1$, we have $c(\mu^{*n}) \le 1$. Furthermore, we can see that $c(\mu^{* (m+n)}) \le c(\mu^{* m}) c(\mu^{*n})$ for any $m,n >0$. Hence, the sequence $ ( c(\mu^{*n}))_{n \ge 1}$ is submultiplicative and satisfies $ c(\mu ^{\ast n_0}) <1$ for some $n_0 \ge 1$.
 
  \noindent Besides, we obtain $m(P^n \varphi) \le m(\varphi) c(\mu^{\ast n})$.  Moreover, we also obtain $m(\varphi) \le |\varphi|_L \le 3m(\varphi)$ for any $\varphi \in Ker \Pi$.  Notice that $P^n(\varphi - \Pi \varphi)$ belongs to $Ker \Pi$ for any $\varphi \in L$ and $n \ge 0$. Hence $|P^n (\varphi - \Pi \varphi)|_L \le 3 c(\mu ^{\ast n}) |\varphi|_L$ which yields
  $$
  |R^n|_{L \to L} = |P^n- \Pi|_{L \to L}  = |P^n(I - \Pi)|_{L \to L} \le 3 c(\mu ^{\ast n}).
  $$
 Therefore, the spectral radius of $R$ is less than or equal to $ \displaystyle \kappa := \lim_{n \to +\infty} \Bigl(c(\mu^{\ast n}) \Bigr)^{1\over n}$ which is strictly less than $1$ by hypothesis P3 and Proposition \ref{propH} (3).
 
   \vspace{0.3cm}
 {\bf (c)} The theory of the perturbation implies that for $\epsilon$ small enough and for any $t \in [-\epsilon; \epsilon]$, the operator $P_t$ may be decomposed as $P _t = \lambda _t \Pi _t + R_t$, where $\lambda _t$ is the dominant eigenvalue of  $P _t$ and the spectral radius of $R_t$ is less than or equal to some $r_\epsilon \in [0,1[$. 
 In order to control $P _t ^n$, we ask $\lambda_t ^n$ to be bounded.  Notice that by Hypothesis P1, the function $t \mapsto P _t$ is analytic near $0$. To prove that the sequence $(P^n_t)_t$ is bounded in $L$, it suffices to check  $|\lambda_t| \le 1$ for any $t \in [-\epsilon, \epsilon]$.
 
 \noindent When $\varphi(x) = {\bf 1} (x)$, equality (\ref{eqn8}) becomes
 \begin{eqnarray} \label{eqn113}
 P _t ^n {\bf 1} (x) = \mathbb E \left[ e^{it\rho( L_n , x)}\right] =  \lambda _t ^n \Pi _t {\bf 1} (x) + R_t ^n {\bf 1} (x) .
 \end{eqnarray}
 
 \noindent We have the local expansion of $\lambda_t$ at $0$: 
 \begin{eqnarray} \label{eqn115}
 \lambda_t = \lambda_0 + t \lambda_0 ' + \frac{t^2}{2} \lambda''_0 [1+ o(1)].
 \end{eqnarray}

 \noindent Taking the first derivative of (\ref{eqn113}) with respect to $t$, we may write for any $n \ge 0$,
 $$
 \frac{d}{dt}  P _t ^n {\bf 1} (x) = \frac{d}{dt} \Bigl( \lambda _t ^n \Pi _t {\bf 1} (x) + R_t^n {\bf 1} (x) \Bigr) = \mathbb E \left[i \rho(L_n , x) e^{it \rho(L_n,x)} \right].
 $$
 Since $\lambda_0 = 1$, $\Pi _0 1(x) =1$ and  $| R^n | _{L \to L} \le C r_\epsilon^n$, we can imply that
 $$
 \lambda_0' =\frac{i}{n} \mathbb E [\rho(L_n,x)] - \frac{\Pi_0' 1(x)}{n} - \frac{[R^n_t 1(x)]'_{t=0}}{n},
 $$
 which yields $  \displaystyle  \lambda_0' = i \lim_{n \to +\infty} {1\over n}\mathbb E [\rho(L_n,x)] = i \gamma _\mu =0 $.  Similarly, taking the second derivative of (\ref{eqn113})
 implies $ \displaystyle \lambda_0'' = - \lim_{n \to +\infty}  \frac{1}{n} \mathbb E [\rho(L_n,x) ^2]$.
 Denote $\sigma ^2 := \displaystyle \lim_{n \to +\infty}  \frac{1}{n} \mathbb E _x [S_n^2]$ . Applying in our context of matrices with non-negative coefficients the argument developed in \cite{BL} Lemma 5.3, we can imply that  $\sigma^2 >0$ and hence $\lambda_0''  = -\sigma^2 <0$. Therefore, in particular, for $t$ closed to $0$, expression (\ref{eqn115}) becomes 
 $$
 \lambda_t= 1 -{\sigma^2\over 2} t^2 [1+ o(1)]
 $$

 which implies $\vert \lambda_t\vert \leq 1$ for $t$ small enough.
 
  \vspace{0.3cm}
 {\bf (d)}  In particular,  inequality  (\ref{controlnormgx}) implies  $\vert \rho(g, x)\vert    \leq \log N(g)$ for any $x \in \mathbb X$. Therefore, for any $p \geq 1$, $x\in \mathbb X$ and  $n\geq 1$, Hypothesis P1 yields
 $$
 \mathbb E_x \vert \rho(g_{n+1}, X_{n})\vert^p  \leq   {p!\over \delta_0^p } \mathbb E_x  e^{\delta_0 \vert \rho(g_{n+1}, X_{n})\vert}  \leq {p! \over \delta_0^p }
 \mathbb E  N(g_{n+1})^{\delta_0}  <+\infty.
 $$
 
 {\rightline{$\square$}}
 
 \subsection{Martingale approximation of the chain $( S_n)_{n \ge 0}$}
 
 \hspace{0.5cm} As announced in the abstract, we approximate the process $(S_n)_{n \ge 0}$ by  a martingale $(M_n)_{n \ge 0}$.  In order to construct the suitable martingale, we introduce the operator $\overline P$ and then find the solution of the Poisson equation as follows.  First, it is neccessary to introduce some notation and  basic properties. Let $g_0 = I$ and $X_{-1} :=X_0 $.
The sequence $((g_n, X_{n-1}))_{n \geq 0}$   is a Markov chain  on $S\times \mathbb X$, starting from $(Id, x)$ and  with   transition operator $\overline P$ defined by: for any $(g, x) \in S \times \mathbb X$ and any bounded measurable function $\phi : S \times \mathbb X \to \mathbb R$,
\begin{eqnarray}\label{eqn27}
\overline {P} \phi(g,x) := \int_{S \times \mathbb X}\phi(h,y) \overline {P} ((g,x), dh dy) = \int_{S } \phi(h,g \cdot x) \mu(dh).
\end{eqnarray}
 
 \noindent 	Notice that, under assumption {\bf P1}, the quantity $\overline P\rho(g, x)$ is well defined for any $(g, x) \in S \times \mathbb X$.

 	 \begin{lem} \label{remark7.3} The function $\displaystyle \bar \rho : x \mapsto \int_S \rho(g, x) \mu(dg) $ belongs to $L$ and for any $g \in S$, $x \in \mathbb X$ and $n \ge 1$,
 	 	\begin{eqnarray} \label{import}
 	 	\overline {P}^{n+1} \rho (g,x) =  P ^n \overline{\rho} (g \cdot x).
 	 	\end{eqnarray}
 	 \end{lem}
 	 {\bf Proof.}
 	 {\bf  (1) }  For any $x \in \mathbb X$, definition of $\rho$ and (\ref{lem2.21}) yield
 	 \begin{eqnarray*}
 	 	|{\overline \rho} (x)| &\leq & \int_{S }  |\log |gx||\mu(dg) \le  \int_{S } 2 |\log N(g)|\mu(dg)  \le  \int_{S } 2 | N(g)|^{\delta_0}\mu(dg) < +\infty. 
 	 \end{eqnarray*}
 	 Hence $|\overline \rho|_\infty < +\infty$. For any $x,y \in \mathbb X$ such that $d(x,y) > \frac{1}{2}$, we can see that 
 	 \begin{eqnarray} \label{eqn821}
 	 |\rho(g,x) - \rho(g,y)| \le |\rho(g,x) - \rho(g,y)| 2 d(x,y) \le 8 \log N(g) d(x,y). 
 	 \end{eqnarray}
 	 For any $x,y \in \mathbb X$ such that $d(x,y) \le \frac{1}{2}$, applying Lemma III.1 in \cite{HH}, we obtain 
 	 \begin{eqnarray} \label{eqn822}
 	 |\rho(g,x) - \rho(g,y)| \le 2 \log \frac{1}{1-d(x,y)} \le 2C d(x,y),
 	 \end{eqnarray}
 	 where $C$ is given in Lemma \ref{lemHH}. For any $x,y \in \mathbb X$, by (\ref{eqn821}) and (\ref{eqn822}) we obtain
 	 \begin{eqnarray*}
 	 	| {\overline \rho} (x) -{\overline \rho} (y)| &\le&  \int_{S} \left|\rho(g,x) - \rho(g,y) \right| \mu(dg) \\
 	 	&\le& \int_{S} [8 \log N(g) +2C ] d(x,y) \mu(dg).
 	 \end{eqnarray*}
 	 Thus $\displaystyle m(\overline \rho) = \sup_{x,y \in \mathbb X, x \ne y} \frac{| {\overline \rho} (x) -{\overline \rho} (y)|}{d(x, y)}  <  +\infty$.
 	 
 	 {\bf  (2)}
 	 From (\ref{eqn27}) and definition of $\rho$, it is obvious that 
 	 \begin{eqnarray*} 
 	 \overline {P}  \rho (g,x) = \int_{S} \rho(h, g \cdot x) \mu(dh)= \overline{\rho} (g \cdot x),
 	 \end{eqnarray*}
 	which yields
 	 \begin{eqnarray*}
 	 	\overline {P}^2  \rho (g,x)&=& \overline {P} (\overline {P} \rho ) (g,x)  = \int_{S \times \mathbb X} (\overline {P} \rho) (k,y) \overline {P}((g,x), dk dy) \\
 	 	&=& \int_{S \times \mathbb X} \overline{\rho}(k \cdot y) \overline {P} ((g,x), dk dy)  \\
 	 	&=&\int_{S} \overline{\rho}(k \cdot (g \cdot x)) \mu(dk) =P \overline{\rho}(g \cdot x).
 	 \end{eqnarray*}
 	 By induction, we obtain $\overline {P}^{n+1} \rho (g,x) =  P ^n \overline{\rho} (g \cdot x)$ for any $n \ge 0$.
 	 
 	 \rightline{$\square$}

 	Formally, the solution $\theta: S\times \mathbb X \to  \mathbb R$ of the equation 
 	$\theta-\overline P \theta = \rho$ is the function 
 	$$
 	\theta : (g, x) \mapsto \sum_{n = 0}^{+\infty} \overline P ^n \rho(g , x).
 	$$
 	Notice that we do not have any spectral property for $\overline P$ and $\rho$ does not belong to $L$. However, we still obtain the convergence of this series by taking into account the important relation  (\ref{import}), as shown in the following lemma.

 	\begin{lem}
 		The sum $\displaystyle \theta =  \sum_{n=0}^{+\infty}  \overline P ^n  \rho$  exists and satisfies the Poisson equation $\rho = \theta - \overline {P} \theta$. Moreover,  		
 	
 			\begin{eqnarray}  \label{lem4.1} 
 			|\overline P \theta |_\infty = \sup_{ g \in S, x \in \mathbb X} |\theta(g,x) - \rho(g,x)| < +\infty;
 			\end{eqnarray}
 		and for any $p \ge 1$, it holds 
 		\begin{eqnarray} \label{remark7.2}
 		\displaystyle \sup_{n \geq 0 } \sup_{x \in \mathbb X} \mathbb E _x| \theta (g_{n+1},X_n)| ^p < +\infty.
 		\end{eqnarray}
 		
 	\end{lem}

 	\noindent{ \bf Proof.} 		
 	{\bf (1)} Since $P$ acts on $(L, |\cdot|_L)$ and $\overline \rho \in L$ from Lemma  \ref{remark7.3}, we obtain $P \overline \rho \in L$. Thanks to  definition of $\rho$, (\ref{eqn10}) and P4, it follows that
 	$$
 	\nu(\overline \rho) = \int_{\mathbb X} \overline \rho (x) \nu(dx) = \int_{S}   \int_{\mathbb X} \rho (g,x) \nu(dx)  \mu(dg) = \gamma _\mu =0.
 	$$
 	 	Proposition \ref{spectre} and the relation (\ref{import}) yield for any  $x \in \mathbb X$ and $n \ge 0$,
 	$$
 	\overline {P} ^{n+1} \rho (g,x) =  P ^n \overline \rho (g \cdot x) = \Pi \overline \rho (g \cdot x) + R^n \overline \rho (g \cdot x)  = \nu (\overline \rho) {\bf 1} (g \cdot x) +R^n \overline \rho (g \cdot x) =  R^n \overline \rho (g \cdot x)
 	$$
 	and there exist $C>0$ and  $0 < \kappa < 1$ such that for any $x \in \mathbb X$ and $n \geq 0$,
 	$$
  \left| R ^n \overline \rho (x) \right| \leq \left| R ^n \overline \rho \right|_L \leq \left| R^n \right|_{L \to L} \leq C \kappa ^n.
 	$$
 	Hence for any $g \in S$ and $x \in \mathbb X$,
 	\begin{eqnarray*}
 		\left| \sum_{n=1}^{+\infty}  \overline P ^n {\rho} (g , x)\right|
 		\le \sum_{n=0}^{+\infty} \left|  P ^n \overline{\rho} (g \cdot x) \right|   \leq C \sum_{n=0}^{+\infty}  \kappa^n = \frac{C}{1-\kappa} < +\infty.
 	\end{eqnarray*}
 	Therefore, the function $\displaystyle \theta =  \sum_{n=0}^{+\infty}  \overline P ^n  \rho$ exists and obviously  satisfies the Poisson equation $\rho = \theta - \overline {P} \theta$.
 	Finally, it is evident that
 	$$
 	\displaystyle \sup_{g \in S, x \in \mathbb X} |\theta(g,x) - \rho(g,x)|=  \sup_{g \in S, x \in \mathbb X} \left| \sum_{n=1}^{+\infty}  \overline P ^n {\rho} (g , x)\right| < +\infty.
 	$$

 {\bf (2)}
 	Indeed, from (\ref{momentsoforderp}), (\ref{lem4.1}) and Minkowski's inequality, the assertion arrives.  
 
 	\rightline{$ \square$}
 	
 
 \vspace{0.5cm}
 Now we contruct a martingale to approximate the Markov walk $(S_n)_{n \geq 0}$. Hence, from the definition of $S_n$  and the Poisson equation, by adding and removing the term $\overline P \theta (g_0, X_{-1})$, we obtain
 \begin{eqnarray*}
 	S_n &=& a+  \rho(g_1, X_0) + \ldots + \rho(g_n,X_{n-1}) \\
 	&=& a+ \overline{P} \theta(g_0, X_{-1}) - \overline{P} \theta (g_n,X_{n-1}) + \sum_{k=0}^{n-1} \left[ \theta (g_{k+1}, X_k) - \overline {P} \theta (g_k, X_{k-1}) \right].
 \end{eqnarray*}
 
 Let $\mathcal F _0:= \{\emptyset, \Omega \}$  and $\mathcal F _n := \sigma \{ g_k: 0\leq k \leq n\}$ for $n \geq1$.

 \begin{prop}
 	For any $n \ge 0$, $x \in \mathbb X$, $a \ge 0$ and $p > 2$, the sequence $(M_n)_{n \geq0}$  defined by 
 	\begin{eqnarray} \label{eqn13}
 	M_0 := S_0\,\, \mbox{and} \,\, M_n:= M_0 + \sum_{k=0}^{n-1} \left[ \theta (g_{k+1}, X_k) - \overline {P} \theta (g_k, X_{k-1}) \right]
 	\end{eqnarray}
 	 is a martingale in $L^p (\Omega, \mathbb P_{x,a}, (\mathcal F _n)_{n\ge 0})$ satisfying the properties:
 	  	 	\begin{eqnarray} \label{lem4.3}
 	 	\displaystyle \sup_{n \geq0} |S_n -M_n| \leq  2 |\overline {P} \theta|_\infty \quad \mathbb P _{x,a}\mbox{-a.s.}
 	 	\end{eqnarray}
 	 
 	  	 	\begin{eqnarray} \label{lem4.4}
 	 	\displaystyle \sup_{n \geq1} n ^{-\frac{p}{2}} \sup_{x \in \mathbb X} \mathbb E _{x,a} |M_n|^{p}< +\infty. 
 	 	\end{eqnarray} 
 
 \end{prop}
 
 { \bf From now on, we set $A:= 2 |\overline {P} \theta|_\infty $.}
 
 \noindent{ \bf Proof.}
 By definition (\ref{eqn13}), martingale property arrives.
 
{\bf (1)} From the construction of $M_n$ and  (\ref{lem4.1}), we can see easily that 
 $$
 \sup_{n \geq0} |S_n - M_n| = \sup_{n \geq0} \left| \overline{P} \theta(g_0, X_{-1}) - \overline{P} \theta (g_n,X_{n-1}) \right| \leq 2 \left| \overline{P} \theta  \right| _\infty < +\infty \quad \mathbb P_{x,a}\mbox{-a.s.}.
 $$

{\bf (2)}
 Denote $\xi_k :=   \theta (g_{k+1}, X_k) - \overline {P} \theta (g_k, X_{k-1}) $. Thus $M_n = M_0 +\sum_{k=0}^{n-1} \xi _k $. 
Using Burkholder's inequality, for any $p \ge 1$, there exists some positive constant $c_p$ such that for $0 \le k <n$,
 $$
 (\mathbb E _{x,a} |M_n |^p)^{\frac{1}{p}} \leq c_p \left( \mathbb E _{x,a} \left| \sum_{k=0}^{n-1} \xi _k^2  \right| ^{\frac{p}{2}} \right) ^{\frac{1}{p}}.
 $$
 Now, with $p >2$, applying Holder's inequality,  we obtain
 
 $$
 \left| \sum_{k=0}^{n-1} \xi _k^2  \right| \le n^{1-\frac{2}{p}} \left(\sum_{k=0}^{n-1} |\xi _k|^p \right)^{\frac{2}{p}},
 $$
 which implies
 
 $$
 \mathbb E _{x,a} \left| \sum_{k=0}^{n-1} \xi _k^2  \right| ^{\frac{p}{2}} \leq n^{\frac{p}{2} -1} \mathbb E _{x,a} \sum_{k=0}^{n-1} |\xi _k|^p \leq n^{\frac{p}{2}} \sup_{0 \leq k \leq n-1} \mathbb E _{x,a} |\xi _k|^p.
 $$
 Since $(M_n)_n$ is a martingale, by using the convexity property, we can see that for any $k \ge 0$,
 $$
 \Big\vert  \overline{P} \theta (g_k, X _{k-1}) \Big\vert ^p = \Big\vert \mathbb E _{x,a}  \Bigl[ |\theta (g_{k+1}, X_k)| |\mathcal F_k \Bigr] \Big\vert ^p \le \mathbb E _{x,a} \Bigl[ |\theta (g_{k+1}, X_k)|^p |\mathcal F_k \Bigr],
 $$
 which implies $ \mathbb E _{x,a} \left| \overline{P} \theta (g_k, X _{k-1})\right| ^p \le \mathbb E _{x,a} \left| \theta(g_{k+1},X_k)\right| ^p  $. 
 Therefore, we obtain 
 \begin{eqnarray*}
 	\Bigl(\mathbb E _{x,a} |M_n |^p\Bigr)^{\frac{1}{p}} &\leq&  c_p\left( n^{\frac{p}{2}} \sup_{0 \leq k \leq n-1} \mathbb E _{x,a} |\xi _k|^p \right) ^{\frac{1}{p}} \le c_p n^{\frac{1}{2}} \sup_{0 \leq k \leq n-1} \Bigl( \mathbb E _{x,a} |\xi_k|^p \Bigr) ^{\frac{1}{p}} \\
 	 &\le&   c_p n^{\frac{1}{2}} \sup_{0 \leq k \leq n-1} \Bigl[ \Bigl( \mathbb E _{x,a} |\theta (g_{k+1}, X_k) |^p \Bigr)^{1/p} + \Bigl( \mathbb E _{x,a} |\overline P \theta (g_k, X_{k-1})|^p \Bigr)^{1/p} \Bigr] \\
 	&\leq&  2 c_p  n^{\frac{1}{2}} \sup_{0 \leq k \leq n-1}  \Bigl( \mathbb E _{x,a} |\theta(g_{k+1}, X_k)|^p \Bigr) ^{\frac{1}{p}}.
 \end{eqnarray*}
 Consequently, we obtain  $\displaystyle \mathbb E _{x,a} |M_n |^p \leq  (2c_p)^p  n^{\frac{p}{2}} \sup_{0 \leq k \leq n-1} \mathbb E _{x,a} |\theta(g_{k+1}, X_k)|^p$ and the assertion arrives by using (\ref{remark7.2}).
 
 \rightline{$ \square$}

\section{On the $\widetilde P _+$-harmonic function $V$ and the proof of Proposition \ref{theo2}}

In this section we construct explicitly a $\widetilde P_+$-harmonic function $V$ and study its properties.  We begin with the first time the martingale $(M_n)_{n \geq 0}$ (\ref{eqn13}) visit $]-\infty, 0]$, defined by 
$$
T = \min \{ n \geq 1: M_n \le 0\}.
$$

The equality $\gamma_\mu=0 $ yields $\displaystyle \liminf_{n \to +\infty} S_n=-\infty \,\, \mathbb P_x$-a.s. for any $x \in \mathbb X$, thus   $\displaystyle \liminf_{n \to +\infty} M_n=-\infty \,\, \mathbb P_x$-a.s., so that   $T <+\infty \,\, \mathbb P_x$-a.s. for any $x \in \mathbb X$ and $ a \ge 0$.

\subsection{On the properties of  $T$ and $(M_n)_n$}

We need to control the first moment  of the random variable $ \vert M_{T \wedge n}\vert$ under $\mathbb P_x$; we consider the restriction of this variable to the event $[T \leq n]$ in lemma  \ref{lem5.2} and  control the remaining term in lemma  \ref{lem5.6}.

\begin{lem}\label{lem5.2}
	
	There exists $\varepsilon _0 >0$ and $c>0$ such that for any $\varepsilon \in (0, \varepsilon _0), n \geq1, x \in \mathbb X$ and $a \geq n^{\frac{1}{2} -\varepsilon}$,
	$$
	\mathbb E _{x,a} \Bigl[|M_{T}|; T \leq n\Bigr] \leq c \frac{a}{n^\varepsilon}.
	$$
\end{lem}
\noindent{ \bf Proof.} 
For any $\varepsilon >0$, consider the event $\displaystyle A_n := \{ \max_{0 \leq k \leq n-1} |\xi_k| \leq n^{\frac{1}{2} -2 \varepsilon} \}$, where $ \xi_k = \theta(g_{k+1},X_k) - \overline{P} \theta(g_k, X_{k-1})$; then
\begin{eqnarray}\label{lelnvkj}
\mathbb E_{x,a} \Bigl[|M_{T}|;T \leq n\Bigr] &=& \mathbb E_{x,a} \Bigl[|M_{T}|;T \leq n,A_n\Bigr] +\mathbb E_{x,a} \Bigl[|M_{T}|;T \leq n,A^c_n\Bigr].
\end{eqnarray}

\noindent On the event $ [T \le n]\cap A_n $, we have  $| M_T | \le |\xi _{T-1}| \le n^{\frac{1}{2} - 2 \varepsilon}$. Hence for any $x \in \mathbb X$ and $a \geq n^{\frac{1}{2} -\varepsilon}$,
\begin{eqnarray} \label{eqn24}
\mathbb E_{x,a} \Bigl[|M_{T}|;T \leq n,A_n\Bigr] \le \mathbb E_{x,a} \Bigl[|\xi_{T-1}|;T \leq n,A_n\Bigr]  \leq n^{\frac{1}{2} -2 \varepsilon} \leq \frac{a}{n^\varepsilon}.
\end{eqnarray}

\noindent Let $M^*_n := \displaystyle \max_{1 \leq k \leq n} |M_k|$; since $|M_{T}| \leq M^*_n$ on the event $[T \leq n]$, it is clear that, for any $x \in \mathbb X$ and $a \ge 0$, 
\begin{eqnarray}\label{Mn*ANc}
	\mathbb E_{x,a} \Bigl[|M_{T}|;T \leq n,A^c_n\Bigr] &\le&  \mathbb E_x [M^*_n; A^c_n] \notag \\
	 &\le& \mathbb E_{x,a} \Bigl[M^*_n;M^*_n > n^{\frac{1}{2} +2\varepsilon}, A^c_n\Bigr] + n^{\frac{1}{2} +2\varepsilon} \mathbb P_{x,a} ( A^c_n) \notag\\
&\le& \int^{+\infty}_{n^{\frac{1}{2}+2\varepsilon}} \mathbb P_{x,a} (M^*_n >t) dt + 2 n^{\frac{1}{2} +2\varepsilon} \mathbb P_{x,a} ( A^c_n).
\end{eqnarray}

\noindent We bound the probability $\mathbb P_{x,a} (A^c_n)$ by using Markov's inequality, martingale definition and (\ref{remark7.2}) as follows: for any $p \ge 1$,
\begin{eqnarray*}
\mathbb P_{x,a} (A^c_n) 
&\le& \sum_{k=0}^{n-1} \mathbb P_{x,a} \left( |\xi_k| > n^{\frac{1}{2} - 2 \varepsilon } \right)\notag \\
&\le& \frac{1}{n^{(\frac{1}{2} - 2 \varepsilon )p}} \sum_{k=0}^{n-1} \mathbb E_{x,a}  |\xi_k|^p \notag \\
&\le& \frac{2^p}{n^{(\frac{1}{2} - 2 \varepsilon )p}} \sum_{k=0}^{n-1} \mathbb E _{x,a} |\theta(g_{k+1}, X_k) |^p      \notag\\
&=& \frac{c_p}{n^{{p \over 2} - 1- 2 \varepsilon p}}.
\end{eqnarray*}
For any $a \ge n^{\frac{1}{2} - \varepsilon}$, it follows that
\begin{equation}\label{TYF}
n^{\frac{1}{2} +2\varepsilon} \mathbb P_{x,a} ( A^c_n) \le  a n^{3\varepsilon}\mathbb P_{x,a} ( A^c_n) \leq \frac{c_p a}{n^{{p \over 2}  -1-2\varepsilon p -3\varepsilon}}.
\end{equation}

Now we control the integral in (\ref{Mn*ANc}). Using Doob's maximal inequality for martingales and (\ref{lem4.4}), we receive for any $p \ge 1$,
	$$
	\mathbb P_x (M^*_n >t) \leq \frac{1}{t^p} \mathbb E_x \Bigl[|M_n|^p \Bigr] \leq c_p \frac{n^{{p \over 2}}}{t^p},
	$$
which implies for any $a \ge n^{\frac{1}{2} - \varepsilon}$,
\begin{eqnarray}\label{eqn5.8}
\int^{+\infty}_{n^{\frac{1}{2}+2\varepsilon}} \mathbb P_x (M^*_n >t) dt 
&\le& {c_p \over p-1} \frac{n^{{p \over 2}}}{n^{(\frac{1}{2}+2\varepsilon)(p-1)}} \le {c_p \over p-1} \frac{a}{n^{2\varepsilon p-3 \varepsilon}}.
\end{eqnarray}
Taking (\ref{Mn*ANc}),    (\ref{TYF})  and (\ref{eqn5.8}) altogether, we obtain for some $c_p'$,
\begin{equation}\label{kzrhd}
\mathbb E_{x,a} \Bigl[|M_{T}|;T \leq n,A^c_n\Bigr] \leq c_p' \left( \frac{a}{n^{2\varepsilon p-3 \varepsilon}} + \frac{ a}{n^{{p \over 2} -1-2\varepsilon p -3\varepsilon}} \right).
\end{equation}

\noindent Finally, from (\ref{lelnvkj}), (\ref{eqn24}) and  (\ref{kzrhd}), we obtain for any $a \geq n^{\frac{1}{2} -\varepsilon}$,
$$
\mathbb E_{x,a} \Bigl[|M_{T}|;T \leq n\Bigr] \leq \frac{a}{n^\varepsilon} +  c_p' \frac{a }{n^\varepsilon}  \left(\frac{1}{n^{2\varepsilon p-4 \varepsilon}} + \frac{ 1}{n^{{p \over 2} -1-2\varepsilon p -4\varepsilon}} \right).
$$
Fix $p >2$. Then there exist $c >0$ and $\varepsilon_0 >0$ such that for any $\varepsilon \in (0, \varepsilon_0) $ and $a \geq n^{\frac{1}{2} -\varepsilon}$,
$$
\mathbb E_{x,a} \Bigl[|M_{T}|, T \le n \Bigr] \leq c \frac{a}{n^\varepsilon}
$$
which proves the lemma.

\rightline{$ \square$}

 For  fixed $\varepsilon >0$ and $a \ge 0$, we consider the first time $\nu_{n, \varepsilon}$ when the process $(|M_k|)_{k \ge 1}$ exceeds $2 n ^{\frac{1}{2} - \varepsilon}$. It is connected to Lemma \ref{lem6.2} where $\mathbb P (\tau_a ^{bm} >n)$ is controlled  uniformly in $a$ under condition $a \le \theta_n \sqrt n$ with $\lim_{ n \to +\infty} \theta_n =0$ which we take into account here by setting
\begin{eqnarray*}
\nu_{n,  \varepsilon} := \min \{k \geq1: | M_k| \geq2 n ^{\frac{1}{2} - \varepsilon} \}.
\end{eqnarray*}
Notice first that  for any $\varepsilon>0, x \in \mathbb X$ and $a \ge 0$ the sequence $( \nu _{n,  \varepsilon})_{n \geq 1}$ tends to $+\infty$ a.s. on $(\Omega, \mathcal B(\Omega), \mathbb P_{x,a})$. The following lemma yields to a more precise control of this property.

\begin{lem}\label{lem5.4}

	For any $\varepsilon \in (0, \frac{1}{2})$, there exists $c_\varepsilon >0$ such that for any $x \in \mathbb X$, $a \ge 0$ and  $n \geq1$,
	$$
	\mathbb P _{x,a} (\nu_{n,  \varepsilon} > n^{1-\varepsilon}) \leq \exp (-c_\varepsilon n^ \varepsilon).
	$$
\end{lem}
\noindent{ \bf Proof.}
Let $m = [B^2 n^{1-2 \varepsilon}]$ and $K = [n^\varepsilon / B^2]$ for some positive constant $B$. By (\ref{lem4.3}), for $n$ sufficiently great  such that $A \leq n^{\frac{1}{2} -\varepsilon}$, we obtain for any $x \in \mathbb X$ and $a \ge 0$,
\begin{eqnarray}\label{eqn5.12}
\mathbb P_{x,a} (\nu_{n,  \varepsilon} > n^{1-\varepsilon}) &\le& \mathbb P_{x,a} \left( \max_{1 \leq k \leq n^{1-\varepsilon}} |M_k| \leq 2n^{\frac{1}{2} -\varepsilon}\right) \notag \\
&\le& \mathbb P_{x,a} \left( \max_{1 \leq k \leq K} |M_{km}| \leq 2n^{\frac{1}{2} -\varepsilon}\right) \notag \\
&\le& \mathbb P_{x,a} \left( \max_{1 \leq k \leq K} |S_{km}| \leq 3n^{\frac{1}{2} -\varepsilon}\right).
\end{eqnarray}
Using  Markov property, it follows that, for any $x \in \mathbb X$ and $a \ge 0$, from which by iterating $K$ times, we obtain
\begin{equation}\label{KJyazdtfhg}
\mathbb P_{x,a} \left( \max_{1 \leq k \leq K} | S_{km}| \leq 3n^{\frac{1}{2} -\varepsilon}\right) \leq \left( \sup_{b \in \mathbb R, x \in \mathbb X} \mathbb P_{x,b} \left(  |S_{m}| \leq 3n^{\frac{1}{2} -\varepsilon}\right) \right) ^K. 
\end{equation}
Denote $\mathbb B (b;r) = \{c: |b+c| \leq r \}$. Then for any $x \in \mathbb X$ and $b \in \mathbb R$
$$
\mathbb P_{x,b} \left( |S_m| \leq 3n^{\frac{1}{2} - \varepsilon }\right) = \mathbb P_x \left( \frac{S_m}{ \sqrt m} \in \mathbb B (b / \sqrt m; r_n) \right),
$$
where $r_n  = {3 n^{\frac{1}{2} -\varepsilon} \over {\sqrt m}}$.  Using the central limit theorem for $S_n$ (Theorem 5.1 property iii) \cite{BL}), we obtain for $n \to +\infty$,
\begin{eqnarray*}
	\sup_{b \in \mathbb R, x \in \mathbb X} \left| \mathbb P_x \left( \frac{S_m}{ \sqrt m} \in \mathbb B  (b/ \sqrt m; r_n) \right) - \int_{\mathbb B (b/ \sqrt m; r_n)} \phi_{\sigma ^2} (u) du \right| \to 0,
\end{eqnarray*}
 where $\phi_{\sigma ^2}(t) = \frac{1}{\sigma\sqrt{2 \pi}} \exp \left( - \frac{t^2}{2 \sigma ^2}\right)$ is the normal density of mean $0$ and variance $\sigma ^2$ on $\mathbb R$.
Since $r_n \leq c_1 B^{-1}$ for some constant $c_1 >0$, we obtain

	$$
	\sup_{b \in \mathbb R} \int_{\mathbb B (b / \sqrt m; r_n)} \phi_{\sigma ^2} (u) du  \leq \int_{-r_n}^{r_n} \phi_{\sigma ^2} (u) du \leq \frac{2r_n}{\sigma \sqrt{2 \pi}} \le \frac{2c_1}{B \sigma \sqrt{2 \pi}}.
	$$

\noindent Choosing $B$  and $n$ great enough, for some $q_\varepsilon <1$, we obtain
$$
\sup_{b \in \mathbb R, x \in \mathbb X} \mathbb P_{x,b} \left( |S_m| \leq 3n^{\frac{1}{2} - \varepsilon }\right) \leq \sup_{b \in \mathbb R} \int_{\mathbb B (b/ \sqrt m; r_n)} \phi_{\sigma ^2} (u) du +o(1)\leq q_\varepsilon .
$$
	Implementing this bound in (\ref{KJyazdtfhg}) and using (\ref{eqn5.12}), it follows that for some $c_\varepsilon >0$,
$$
\sup_{a >0 , x \in \mathbb X} \mathbb P_{x,a} (\nu_{n, \varepsilon} >n^{1-\varepsilon}) \leq q_\varepsilon ^K  \le q_\varepsilon^{\frac{n^\varepsilon}{B^2} -1}  \leq e^{- c_\varepsilon n^\varepsilon }.
$$


\rightline{$ \square$}

\begin{lem}\label{lem5.5} 
	
	There exists $c >0$ such that for any $\varepsilon \in (0, \frac{1}{2})$, $ x \in \mathbb X$, $a \ge 0$ and  $n \geq1$,
	
		$$
		\sup_{1 \leq k \leq n} \mathbb E _{x,a} [|M_k|; \nu_{n, \varepsilon} > n ^{1- \varepsilon}] \leq c (1+a) \exp (-c _\varepsilon n^ \varepsilon)
		$$
	
	for some positive constant $c_\varepsilon $ which only depends on $\varepsilon$.
\end{lem}
\noindent{ \bf Proof.}
By Cauchy-Schwartz inequality, for any $x \in \mathbb X$, $a \ge 0$ and $1 \leq k \leq n$,
$$
\mathbb E_{x,a} \left[|M_k|; \nu_{n, \varepsilon} > n^{1-\varepsilon}\right] \leq \sqrt{ \mathbb E_{x,a} |M_k|^2 \mathbb P_{x,a} ( \nu_{n, \varepsilon} > n^{1-\varepsilon})}.
$$
By Minkowsky's inequality, (\ref{lem4.3}) and the fact that $\frac{1}{n} \mathbb E_x |S_n|^2 \to \sigma^2$ as $n \to +\infty$, it yields
$$
\sqrt{ \mathbb E_{x,a} |M_k|^2}  \leq a+ \sqrt{  \mathbb E_{x,a} [M_k^2]}  \leq a+ A + \sqrt{ \mathbb E_{x,a} [S_k^2]}  \leq c (a+ n^ {\frac{1}{2}})
$$
for some $c >0$ which does not depend on $x$.	The claim follows by Lemma \ref{lem5.4}.

\rightline{$ \square$}

\begin{lem}\label{lem5.6}
	
	There exists $c >0$ and $\varepsilon_0 >0$ such that for any $\varepsilon \in (0, \varepsilon_0)$,  $x \in \mathbb X$, $a \ge 0$ and $n \geq1$,
	\begin{equation}	\label{5.15}
		\mathbb E _{x,a} [M_n; T >n] \leq c(1+a).
		\end{equation}
		and
		 \begin{equation} \label{lem5.10}
		\lim_{a \to +\infty} \frac{1}{a} \lim_{n \to +\infty} \mathbb E _{x,a} [M_n;T >n]= 1.
		\end{equation}  
			
\end{lem}
\noindent{ \bf Proof.} 
{ \bf (1)} \noindent On one hand, we claim 
\begin{eqnarray} \label{eqnclaim}
\mathbb E_{x,a} [M_n; T >n , \nu_{n,  \varepsilon} \leq n^{1-\varepsilon}]	&\leq &  \left( 1+ \frac{c'_\varepsilon}{n^\varepsilon} \right)  \mathbb E_{x,a} \left[M_{[n^{1-\varepsilon}]}; T >[n^{1-\varepsilon}]\right]
\end{eqnarray}
and delay the proof of (\ref{eqnclaim}) at the end of the first part.
\noindent On the other hand, by Lemma \ref{lem5.5}, there exists $c >0$ such that for any $\varepsilon \in (0, \frac{1}{2})$,  $x \in \mathbb X$, $a \ge 0$ and $n \ge 1$,
\begin{eqnarray}\label{eqn5.112}
\mathbb E_{x,a} [M_n; T >n , \nu_{n,  \varepsilon} > n^{1-\varepsilon}] &\le& \sup_{1 \le k \le n} \mathbb E_{x,a} \Bigl[|M_k|; \nu_{n,  \varepsilon} > n^{1-\varepsilon}\Bigr] \notag \\
&\le&  c(1+a) \exp (-c_\varepsilon n^\varepsilon).
\end{eqnarray}

\noindent Hence  combining (\ref{eqnclaim}) and (\ref{eqn5.112}), we obtain for any $x \in \mathbb X$ and $a \ge 0$,
\begin{equation}\label{5.23}
\mathbb E_{x,a} [M_n; T >n] \leq \left( 1+ \frac{c'_\varepsilon}{n^\varepsilon} \right)  \mathbb E_{x,a} \left[M_{[n^{1-\varepsilon}]}; T >[n^{1-\varepsilon}]\right] +c(1+a) \exp (-c_\varepsilon n^\varepsilon).
\end{equation}

\noindent Let $k_j := \left[ n^{(1-\varepsilon)^j} \right]$ for $j \geq0$. Notice that $k_0 =n$ and $[k_j^{1-\varepsilon}] \le  k_{j+1}$ for any $j \ge 0$. Since 
the sequence  $((M_n){\bf 1}_{[T >n]})_{n\geq1}$ is a submartingale,  by using the bound (\ref{5.23}), it yields
\begin{eqnarray*}
	 \mathbb E_{x,a} [M_{k_1}; T >{k_1}] 
	&\le&  \left( 1+ \frac{c'_\varepsilon}{{k_1}^\varepsilon} \right) \mathbb E_{x,a} \left[M_{[{k_1}^{1-\varepsilon}]}; T >[{k_1}^{1-\varepsilon}]\right] + c(1+a) \exp (-c_\varepsilon k_1^\varepsilon) \\
	&\le& \left( 1+ \frac{c'_\varepsilon}{{k_1}^\varepsilon} \right) \mathbb E_{x,a} [ M_{k_2}; T >k_2] + c(1+a) \exp (-c_\varepsilon k_1^\varepsilon).
\end{eqnarray*}
Let $n_0$ be a constant and $m = m(n)$  such that $k_m = \left[ n^{(1-\varepsilon )^m}\right] \le n_0$. After $m$ iterations, we obtain
\begin{equation}\label{5.24}
\mathbb E_{x,a} [M_n; T >n] \leq A_m \Bigl( \mathbb E_{x,a} [M_{k_m}; T >k_m] + c(1+a) B_m \Bigr),
\end{equation}
where 
	\begin{eqnarray}\label{5.25}
	A_m = \prod _{j=1}^m \left( 1+ \frac{c'_\varepsilon}{k_{j-1}^\varepsilon}\right) \le \exp \Bigl( {2^\varepsilon c'_\varepsilon \frac{n_0^{-\varepsilon}}{ 1- n_0^{-\varepsilon^2} }} \Bigr),
	\end{eqnarray}
and 
	\begin{equation}\label{5.26}
	B_m = \sum_{j=1}^m \frac{\exp \left( -c_\varepsilon k_{j-1}^\varepsilon \right)}{\Bigl( 1+\frac{c'_\varepsilon}{k^\varepsilon_{j-1}} \Bigr) \ldots \Bigl( 1+\frac{c'_\varepsilon}{k^\varepsilon_m} \Bigr)   } \leq c_1 \frac{n_0^{-\varepsilon}} {1- n_0^{-\varepsilon^2}}
	\end{equation}
from Lemma 5.6 in \cite{GLP1}. By choosing $n_0$ sufficient great, the first assertion of the lemma follows  from (\ref{5.24}), (\ref{5.25})
and (\ref{5.26}) taking into account that
	\begin{eqnarray*}\label{5.29}
	\mathbb E_{x,a} [ M_{k_m}; T >k_m] \le  \mathbb E_{x,a} [ M_{n_0}; T > n_0] \leq \mathbb E_{x,a}  |M_{n_0}|  \leq a+c.
	\end{eqnarray*}

\noindent Before proving (\ref{eqnclaim}), we can see that  there exist $c >0$ and $0 < \varepsilon_0 < {1 \over 2}$ such that for any  $\varepsilon \in (0, \varepsilon_0)$, $x \in \mathbb X$ and $b \geq n^{\frac{1}{2} -\varepsilon}$,
\begin{equation}\label{jktbzh}
\mathbb E_{x,b} [M_n;T >n] \leq \left( 1+\frac{c}{n^\varepsilon}\right)b.
\end{equation}
Indeed, since $(M_n, \mathcal F_n)_{n \geq1}$ is a  $\mathbb P_{x,b}$- martingale, we obtain 
$$\mathbb E_{x,b} [M_n;T \leq n] = \mathbb E_{x,b} [M_{T};T \leq n]$$ 
and thus
\begin{eqnarray}\label{jkjejkva}
\mathbb E_{x,b} [M_n;T>n]  &=& \mathbb E_{x,b} [M_n] -\mathbb E_{x,a} [M_n;T \leq n]\notag \\
&=&  b -  \mathbb E_{x,b} [M_{T};T \leq n] \notag\\
&=&  b +  \mathbb E_{x,b} [|M_{T}|;T \leq n] .
\end{eqnarray}
Hence (\ref{jktbzh}) arrives by using Lemma \ref{lem5.2}.
For  (\ref{eqnclaim}), it is obvious that
\begin{equation}\label{5.19}
\mathbb E_{x,a} \Bigl[M_n; T >n , \nu_{n,  \varepsilon} \leq n^{1-\varepsilon}\Bigr]  = \sum_{k=1}^{[n^{1-\varepsilon}]} \mathbb E_{x,a} \Bigl[M_n; T >n, \nu_{n,  \varepsilon} =k\Bigr].
\end{equation}
Denote $U_m(x, a) := \mathbb E_{x,a} [M_m; T >m]$. For  any $m \geq1$, by the Markov property applied to $(X_n)_{n \geq1}$, it follows that
\begin{eqnarray}\label{hhgfjkegr}
\mathbb E_{x,a} \Bigl[M_n; T >n, \nu_{n,  \varepsilon} =k\Bigr] &=& \int \mathbb E_{y,b} [M_{n-k}; T>n-k] \notag\\
&& \qquad \qquad  \mathbb P_{x,a} (X_k \in dy, M_k \in db; T >k, \nu_{n, \varepsilon }=k)\notag \\
&=& \mathbb E_{x,a} \Bigl[U_{n-k} (X_k, M_k); T >k,\nu_{n, \varepsilon } =k \Bigr].
\end{eqnarray}
 From the definition of $\nu_{n, \varepsilon}$, we can see that $[\nu_{n, \varepsilon} =k] \subset \left[ |M_k|\geq n^{\frac{1}{2} - \varepsilon} \right]$, and by using (\ref{jktbzh}), on the event $[T >k, \nu_{n,  \varepsilon} =k]$ we have $U_{n-k} (X_k, M_k) \leq \left( 1 + \frac{c}{(n-k) ^\varepsilon} \right) M_k$.
Therefore (\ref{hhgfjkegr}) becomes
\begin{eqnarray}\label{5.22}
\mathbb E_{x,a} [M_n; T >n, \nu_{n, \varepsilon}=k] \leq \left( 1+ \frac{c}{(n-k)^\varepsilon} \right) \mathbb E_{x,a} [M_k; T >k, \nu_{n,  \varepsilon}=k].
\end{eqnarray}
Combining (\ref{5.19}) and (\ref{5.22}), it follows that, for $n$ sufficiently great,
\begin{eqnarray*}
	\mathbb E_{x,a} [M_n; T >n , \nu_{n,  \varepsilon} \leq n^{1-\varepsilon}] &\leq &  \sum_{k=1}^{[n^{1-\varepsilon}]} \left( 1+ \frac{c}{(n-k)^\varepsilon} \right) \mathbb E_{x,a} [M_k; T >k, \nu_{n,  \varepsilon} =k]\\
	&\le& \left( 1+ \frac{c'_\varepsilon}{n^\varepsilon} \right)  \sum_{k=1}^{[n^{1-\varepsilon}]}  \mathbb E_{x,a} [M_k; T >k, \nu_{n, \varepsilon} =k],
\end{eqnarray*}
for some constant $c'_\varepsilon >0$. Since $(M_n{\bf 1}_{[T >n]})_{n\geq1}$ is a submartingale, for any $x \in \mathbb X$, $a \ge 0$ and $1 \leq k \leq [n^{1-\varepsilon}]$,
$$
\mathbb E_{x,a} [M_k; T >k, \nu_{n,  \varepsilon} =k]  \leq \mathbb E_{x,a} \Bigl[M_{[n^{1-\varepsilon}]}; T >[n^{1-\varepsilon}], \nu_{n, \varepsilon} =k\Bigr].
$$
This implies 
\begin{eqnarray*}
	\mathbb E_{x,a} [M_n; T >n , \nu_{n,  \varepsilon} \leq n^{1-\varepsilon}] &\leq &  \left( 1+ \frac{c'_\varepsilon}{n^\varepsilon} \right) \sum_{k=1}^{[n^{1-\varepsilon}]}   \mathbb E_{x,a} \left[M_{[n^{1-\varepsilon}]}; T >[n^{1-\varepsilon}], \nu_{n, \varepsilon} =k\right]\\
	&\leq &  \left( 1+ \frac{c'_\varepsilon}{n^\varepsilon} \right)  \mathbb E_{x,a} \left[M_{[n^{1-\varepsilon}]}; T >[n^{1-\varepsilon}]\right].
\end{eqnarray*}


{\bf (2)} Let $\delta >0$. From  (\ref{5.25}) and (\ref{5.26}), by choosing $n_0$ sufficiently great, we obtain $A_m \leq 1+\delta$ and $B_m \leq \delta$. Together with (\ref{5.24}), since $(M_n {\bf 1}_{[T >n]})_{n \geq1}$ is a submartingale, we obtain for $k_m \le n_0$,
$$
\mathbb E_{x,a} [M_n; T >n] \leq (1+\delta) \Bigl( \mathbb E_{x,a} [ M_{n_0}; T > n_0] +c(1+a) \delta   \Bigr).
$$
Moreover, the sequence $\mathbb E_{x,a} [M_n; T >n]$ is increasing, thus it  converges $\mathbb P_{x,a}$-a.s. and 
$$
\lim_{n \to +\infty} \mathbb E_{x,a} [M_n; T >n] \leq (1+\delta) \Bigl( \mathbb E_{x,a} [M_{n_0}; T >n_0] +c(1+a) \delta \Bigr).
$$
By using (\ref{jkjejkva}), we obtain
$$
a \leq \lim_{n \to +\infty} \mathbb E_{x,a} [M_n; T >n] \leq (1+ \delta) \Bigl(  a+ \mathbb E_x |M_{n_0}| + c(1+a) \delta \Bigr).
$$ Hence the assertion follows since $\delta >0$ is arbitrary.

\rightline{$ \square$}

\subsection{On the stopping time  $\tau$.}

We now state some useful properties of $\tau$ and $S_\tau$.
\begin{lem}\label{lem5.7}
	
	There exists $c>0$ such that for any  $x \in \mathbb X$, $a \ge 0$ and $n \geq1$,
	\begin{equation*}
	\mathbb E _{x,a} [ S_n , \tau >n  ]  \leq  c(1+a).   
	\end{equation*}
\end{lem}
{\bf Proof.} (\ref{lem4.3}) yields $\mathbb P _{x} (\tau_a \le T_{a+A}) =1$ and  $ A+M_n \geq S_n >0$ on the event $[\tau >n]$. By (\ref{5.15}), it follows that 
\begin{eqnarray*}
\mathbb E_{x,a} [S_n; \tau >n] &\leq& \mathbb E_{x,a} [ A+M_n ; \tau >n] \notag \\  
&\leq&  \mathbb E_{x,a+A} [ M_n; T >n] \\
&\leq& c_1 (1+ a+ A) \leq c_2(1+a).
\end{eqnarray*}

\rightline{$ \square$}


\begin{prop} \label{lem5.9} 
	
	There exists $c >0$ such that for any $x \in \mathbb X$ and $a \ge 0$,
	 \begin{eqnarray*}
			\mathbb E_{x,a} |S_\tau| \le c(1+a) < +\infty,
		\end{eqnarray*}
	 \begin{eqnarray} \label{lem3.72}
		\mathbb E_{x,a} |M_\tau| \le c(1+a) < +\infty.
		\end{eqnarray}

\end{prop}
{\bf Proof.} By (\ref{lem4.3}), since $(M_n)_n$ is a martingale, we can see that 
\begin{eqnarray*}
	-\mathbb E_{x,a} [S_{\tau} ; \tau \leq n] &\le& -\mathbb E_{x,a} [M_{\tau} ; \tau \leq n] +A \\
	&=& \mathbb E_{x,a} [M_n ; \tau >  n] -\mathbb E_{x,a} [M_n] +A \\
	&\le& \mathbb E_{x,a} [S_n ; \tau >  n]  + 2A.
\end{eqnarray*}
Hence by Lemma \ref{lem5.7}, for any $x \in \mathbb X$ and $a \ge 0$,
\begin{eqnarray*}\label{5.31}
\mathbb E_{x,a} \left[ |S_{\tau }| ;\tau \le n \right]  &\leq&  \mathbb E_{x,a} \left| S_{\tau \wedge n} \right| \notag \\
 &=& \mathbb E_{x,a} \left[ S_{n } ;\tau > n \right]  - \mathbb E_{x,a} \left[ S_{\tau } ;\tau \leq n \right] \notag \\
&\le& 2 \mathbb E_{x,a} [S_n ; \tau >  n]  + 2A \notag \\
&\le&   c(1+a) +2A.
\end{eqnarray*}

\noindent  By Lebesgue's Dominated Convergence Theorem, it yields 
$$
\mathbb E_{x,a}  |S_{\tau }| = \lim_{n \to +\infty}   \mathbb E_{x,a} \left[ |S_{\tau}| ;\tau \le n \right]  \leq  c(1+a) +2A < +\infty.
$$
By (\ref{lem4.3}), the second assertion arrives.

\rightline{$ \square$}


\subsection{Proof of Proposition \ref{theo2}}
 Denote $\tau_a := \min \{n \ge 1: S_n \le -a \}$ and $T_a := \min \{n \ge 1: M_n \le -a \}$  for any $a \ge 0$. Then $\mathbb E _{x,a} M _{\tau} = a+ \mathbb E _{x} M _{\tau_a} $ and $\mathbb P_{x,a} (\tau >n) = \mathbb P_{x} (\tau_a > n) $.

{\bf (1)} By (\ref{lem3.72}) and Lebesgue's Dominated Convergence Theorem, for any $x \in \mathbb X$ and $a \ge 0$,
$$
\lim_{n \to +\infty} \mathbb E_{x,a} [ M_{\tau}; \tau \leq n ] = \mathbb E_{x,a}  M_{\tau} =a - V(x,a),
$$
where $V(x,a)$ is the quantity defined  by: for $x \in \mathbb X$ and $a \in \mathbb R$,
\begin{eqnarray*}
	V(x,a) := \left\{ 
	\begin{array}{l}
		- \mathbb E _x M_{\tau _a} \quad \mbox{if} \quad a \ge 0,\\
		0 \qquad \qquad \,\, \mbox{if} \quad a < 0.
	\end{array} \right.
\end{eqnarray*}

 \noindent Since $(M_n, \mathcal F_n)_{n \geq1}$ is a $\mathbb P_{x,a}$-martingale, 
\begin{eqnarray}\label{eqn22}
\mathbb E_{x,a} [M_n; \tau >n] &=& \mathbb E_{x,a} M_n - \mathbb E_{x,a} [M_n; \tau \leq n] = a - \mathbb E_{x,a} [M_\tau; \tau \leq n],
\end{eqnarray}
which implies
$$
\lim_{n \to +\infty} \mathbb E_{x,a} [ M_n; \tau > n ] = V(x, a).
$$
Since $\left| S_n - M_n \right| \leq A$ $\mathbb P _x$-a.s. and $\displaystyle \lim_{n \to +\infty} \mathbb P_{x,a} (\tau >n) =0$, it follows that
$$
\lim_{n \to +\infty} \mathbb E_{x,a} [S_n; \tau >n] = \lim_{n \to +\infty} \mathbb E_{x,a} [M_n; \tau >n]  = V(x, a).
$$


{\bf (2)} The assertion arrives by taking into account that $0 \le a \leq a'$ implies $\tau_a \leq \tau_{a'}$ and 
$$
\mathbb E_x [a+S_n; \tau_a >n] \leq  \mathbb E_x [a' +S_n; \tau_{a'} >n].
$$

{\bf  (3)} Lemma \ref{lem5.7} and assertion 1  imply that $V(x, a) \leq c(1+a)$  for any $x \in \mathbb X$ and $a \ge 0$. Besides, (\ref{eqn22}) and (\ref{lem4.3}) yield
$$
\mathbb E_{x,a} [M_n ; \tau >n] \geq a - \mathbb E_{x,a} [ S_\tau ; \tau \leq n] - A  \geq a -A,
$$
which implies 
\begin{eqnarray} \label{eqn231}
V(x,a) \ge a-A.
\end{eqnarray}
 Now we prove $V(x,a) \ge 0$. Assertion 2 implies $V(x,0) \le V(x,a)$ for any $x \in \mathbb X$ and $a \ge 0$. From P5, let $E_\delta := \{g \in S: \forall x \in \mathbb X, \log|gx| \ge \delta \}$ and choose a positive constant $k$ such that $k\delta > 2A$. Hence, for any $g_1, \ldots, g_k \in E_\delta$ and any $x \in \mathbb X$, we obtain $\log |g_k \ldots g_1 x| \ge k \delta > 2A$. It yields
\begin{eqnarray*}
	V(x,0) &=& \lim_{n \to +\infty} \mathbb E_{x} [S_n ; \tau >n] \\
	&\ge & \liminf_{n \to +\infty} \int_{E_\delta} \ldots \int_{E_\delta} \mathbb E_{g_k \ldots g_1 \cdot x, \log |g_k \ldots g_1x|} [S_{n-k} ; \tau >n-k] \mu(dg_1) \ldots \mu(dg_k) 	\\
	&\ge & \liminf_{n \to +\infty} \int_{E_\delta} \ldots \int_{E_\delta} V(g_k \ldots g_1 \cdot x,2A) \mu(dg_1) \ldots \mu(dg_k) 	\\
	&\ge& A \Bigl( \mu(E_\delta) \Bigr)^k >0,
\end{eqnarray*}
where the last inequality comes from (\ref{eqn231}) by applying to $a = 2A$.

 
 {\bf (4)} Equation (\ref{eqn231}) yields $\displaystyle \lim_{a \to +\infty} \frac{V(x, a)}{a} \geq1$. By (\ref{lem4.3}), it yields that $  \mathbb P_{x} (\tau_a < T_{A+a}) =1$, which implies
$$
\mathbb E_{x,a} [S_n; \tau >n] \leq \mathbb E_{x,a} [A+ M_n ; \tau >n] \leq \mathbb E_{x,a} [ A +M_n; T_{ A} >n] = \mathbb E_{x,a+A} [M_n; T >n].
$$
From (\ref{lem5.10}), we obtain $\displaystyle\lim_{n \to +\infty} \frac{V(x, a)}{a} \leq 1$.


{ \bf (5)}  For any $x \in \mathbb X$, $a \ge 0$ and $n \geq 1$, we set $V_n (x, a) := \mathbb E_{x,a} [S_n; \tau >n]$. By assertion 1, we can see $\displaystyle \lim_{n \to +\infty} V_n(x,a) = V(x,a)$. By Markov property, we obtain
\begin{eqnarray*}
	V_{n+1}(x, a)
	&=& \mathbb E_{x,a} \left[  \mathbb E \left[  S_1 +  \sum_{k=1}^{n} \rho(g_{k+1} ,X_k) ; S_1  >0, \ldots ,  S_{n+1} >0 |\mathcal F_1 \right]    \right] \\
				&=& \mathbb E_{x,a} \left[  V_n(X_1,S_1)  ;\tau >1    \right]. 
\end{eqnarray*}
By Lemma \ref{lem5.7}, we obtain $\displaystyle \sup_{x \in \mathbb X,a \ge 0} V_n(x, a) \leq c(1+a)$ which implies $\mathbb P$-a.s.
$$
V_n(X_1,S_1) {\bf 1}_{[\tau >1]} \leq  c(1+S_1){\bf 1}_{[\tau >1]}. 
$$
 Lebesgue's Dominated Convergence Theorem and (\ref{eqn7}) yield
\begin{eqnarray*} \label{eqn23}
	V(x, a) = \lim_{n \to +\infty} V_{n+1} (x,a) &=&  \lim_{n \to +\infty}  \mathbb E_{x,a} [ V_n (X_1, S_1); \tau >1] \notag \\
	&=& \mathbb E_{x,a} [ V (X_1, S_1); \tau >1] \notag \\
	&=& \widetilde P _+ V(x, a). 
\end{eqnarray*}
 {\rightline{$\square$}}

\section{Coupling argument and proof of theorems \ref{theo3} and \ref{theo4}}

\hspace{0.5cm} First, we apply the weak invariance principle stated in \cite{GLP1} and verify that the sequence $(\rho(g_k, X_{k-1}))_{k \ge 0}$ satisfies the conditions of Theorem 2.1 \cite{GLP1}. The hypotheses C1, C2 and C3 of this theorem are given in terms of Fourier transform of the partial sums of $S_n$. Combining the expressions (\ref{eqn201}), (\ref{eqn8}), (\ref{eqn1.131}) and the properties of the Fourier operators $(P_t)_t$, we verify in the next section that the conditions C1, C2 and C3 of Theorem 2.1 in \cite{GLP1} are satisfied in our context.  This lead to the following simpler statement but sufficient. 
\begin{theo}
	Assume P1-P4. There exist
	\begin{itemize}
		\item $\varepsilon_0 >0$, and $c_0 >0$,
		
		\item a probability space $(\widetilde \Omega, B(\widetilde \Omega))$,
		
		\item a family $(\widetilde {\mathbb P} _x)_{x \in \mathbb X}$ of probability measures on $(\widetilde \Omega, B(\widetilde \Omega))$,
		
		\item a sequence $(\tilde a _k)_k$ of real-valued random variables on $(\widetilde \Omega, B(\widetilde \Omega))$  such that $\mathcal L \Bigl((\tilde a _k)_k \slash \widetilde {\mathbb P} _x\Bigr) = \mathcal L \Bigl(( a _k)_k \slash  {\mathbb P} _x\Bigr) $ for any $x \in \mathbb X$,
		
		\item and a sequence $(\widetilde W _i)_{i \ge 1}$ of independent standard normal random variables on $(\widetilde \Omega, \mathcal B(\widetilde \Omega)) $
		
	\end{itemize}
	
	such that for any $x \in \mathbb X$,
	\begin{equation} \label{eqn212}
	\widetilde{\mathbb{P} _x}\left(\sup_{1 \le k\leq n}\left\vert  	\sum_{i=1}^{k}(\tilde{a}_i -\sigma \widetilde W_i)
	\right\vert >  n^{{1\over 2}-\varepsilon_0}\right) \leq c_0 n^{- \varepsilon_0}.
	\end{equation}

\end{theo}

Notice that the fact (\ref{eqn212}) holds true for $\varepsilon_0$ implies (\ref{eqn212}) holds true for $\varepsilon $, whenever $\varepsilon  \le \varepsilon _0$. In order to simplify the notations, we identify $(\widetilde \Omega, \mathcal B( \tilde \Omega))$ and $( \Omega, \mathcal B(  \Omega))$ and consider that  
 the process $(\log |L_n x|)_{n \ge 0}$  satisfies the following property: there exists $\varepsilon_0 >0$ and $c_0 >0$ such that for any $\varepsilon \in (0, \varepsilon_0]$ and $x \in \mathbb X$,
\begin{eqnarray}\label{eqn6.5}
\mathbb P \left( \sup_{0 \le t \le 1} | \log |L_{[nt]}x| - \sigma B_{nt}| > n^{\frac{1}{2} -\varepsilon}\right) = 
\mathbb P_x \left( \sup_{0 \le t \le 1} | S_{[nt]} - \sigma B_{nt}| > n^{\frac{1}{2} -\varepsilon}\right) \le c_0 n^{-\varepsilon},
\end{eqnarray}
where $(B_t)_{t \ge 0}$ is a standard Brownian motion on the probability space $(\Omega, \mathcal B(\Omega), \mathbb P)$ and $\sigma >0$ is used in the proof of Proposition \ref{spectre}, part c). For any $a \ge 0$, let $\tau_a^{bm}$ be the first time the process $(a+ \sigma B_t)_{t\ge 0}$ becomes non-positive: 
$$
\tau_a^{bm} = \inf \{t \ge 0: a+\sigma B_t \le 0 \}.
$$

  The following lemma is due to Levy \cite{Levy} (Theorem 42.I, pp. 194-195).
  \begin{lem} \label{lem6.1}
  	
  	\begin{enumerate}
  		\item  For any $a\ge 0$ and $n \ge 1$, 
  		\begin{eqnarray*} 
  		\mathbb P (\tau_a^{bm} >n)  = \mathbb P \left( \sigma \inf_{0 \le u \le n} B_u > -a \right) =  \frac{2}{\sigma \sqrt {2 \pi n}} \int_0^a \exp \left( - \frac{s^2}{2n \sigma^2}\right) ds.
  		\end{eqnarray*} 
  		
  		\item  For any $a, b$ such that $0 \le a < b< +\infty$ and $n \ge 1$,
  		\begin{eqnarray*} 
  		&& \mathbb P(\tau_a^{bm} >n, a + \sigma B_n \in [a,b]) \notag \\
  		 && \qquad = \frac{1}{\sigma \sqrt {2 \pi n}} \int_a^b \left[ \exp \left( - \frac{(s-a)^2}{2n \sigma^2}\right) -\exp \left( - \frac{(s+a)^2}{2n \sigma^2}\right) \right] ds.
  		\end{eqnarray*}
  		
  	\end{enumerate}
  \end{lem}
From lemma \ref{lem6.1}, we can obtain the next result.  
  
\begin{lem} \label{lem6.2}

	\begin{enumerate}
		\item  There exists a positive constant $c$ such that for any $a\ge 0$ and $n \ge 1$, 
		\begin{eqnarray} \label{eqn6.3}
		\mathbb P (\tau_a^{bm} >n) \le c\frac{a}{\sigma \sqrt n}.
		\end{eqnarray} 
		
		\item  For any sequence of real numbers $(\alpha_n)_n$ such that $\alpha_n \to 0$ as $n\to +\infty$, there exists a positive constant $c$ such that for any $a \in [0 , \alpha_n \sqrt n]$,
		\begin{eqnarray} \label{eqn6.4}
	 \Big| \mathbb P(\tau_a^{bm} >n) - \frac{2a}{\sigma \sqrt{2 \pi  n}} \Big|  \le  c\frac{\alpha_n}{\sqrt n} a.
		\end{eqnarray}
	
	\end{enumerate}
\end{lem}

We use the coupling result described in Theorem 4.1 above to transfer the properties of the exit time $\tau_a^{bm}$ to the exit time $\tau_a$ for  great $a$.

\subsection{Proof of Theorem \ref{theo3}}

\hspace{0.5cm} {\bf (1)} Let $\varepsilon \in (0, \min\{\varepsilon_0; {1 \over 2}\}) $ and $(\theta_n)_{n \ge 1}$ be a sequence of positive numbers such that $\theta_n \to 0$ and $\theta_n n^{\varepsilon /4} \to +\infty$ as $n \to +\infty$. For any  $x \in \mathbb X$ and $ a  \ge 0$, we have the decomposition
\begin{eqnarray} \label{eqn6.13}
P_n (x,a) := \mathbb P _{x,a} (\tau >n) &=& \mathbb P  _{x,a} (\tau >n, \nu_{n,\varepsilon} > n^{1 - \varepsilon})  +\mathbb P _{x,a} (\tau >n, \nu_{n,\varepsilon} \le n^{1 - \varepsilon}) .
\end{eqnarray}
It is obvious that from lemma \ref{lem5.4}, we obtain
\begin{eqnarray} \label{eqn6.14}
\sup_{x \in \mathbb X, a \ge 0} \mathbb P _{x,a} (\tau >n, \nu_{n,\varepsilon} > n^{1 - \varepsilon} ) \le \sup_{x \in \mathbb X, a \ge 0} \mathbb P _{x,a} ( \nu_{n,\varepsilon} > n^{1 - \varepsilon} )  \le e^{-c_\varepsilon n^\varepsilon}.
\end{eqnarray}
For the second term, by Markov's property,
\begin{eqnarray} \label{eqn6.15}
  \mathbb P _{x,a} (\tau >n, \nu_{n,\varepsilon} \le n^{1 - \varepsilon}) 
 &= &   \mathbb E_{x,a} \left[ P _{n - \nu _n} (X_{\nu_{n,\varepsilon}}, S_{\nu_{n,\varepsilon}}); \tau > \nu_{n,\varepsilon}, \nu_{n,\varepsilon} \le n^{1-\varepsilon} \right] \label{4.67} \\
  & =&  I_n(x,a) + J_n(x,a), \notag 
\end{eqnarray}
where
\begin{eqnarray*}
 I_n(x,a) := \mathbb E_{x,a} \left[ P _{n - \nu _n} (X_{\nu_{n,\varepsilon}}, S_{\nu_{n,\varepsilon}}); S_{\nu_{n,\varepsilon}} \le \theta_n n^{\frac{1}{2}} , \tau > \nu_{n,\varepsilon}, \nu_{n,\varepsilon} \le n^{1-\varepsilon} \right], \\
\mbox{and}\,  \,J_n(x,a) := \mathbb E_{x,a} \left[ P _{n - \nu _n} (X_{\nu_{n,\varepsilon}}, S_{\nu_{n,\varepsilon}}); S_{\nu_{n,\varepsilon}} > \theta_n n^{\frac{1}{2}} , \tau > \nu_{n,\varepsilon}, \nu_{n,\varepsilon} \le n^{1-\varepsilon} \right]. 
\end{eqnarray*}

Now we control the quantity $P _{n - \nu _n} (X_{\nu_{n,\varepsilon}}, S_{\nu_{n,\varepsilon}})$ by using the following lemma. The proofs of the lemmas stated in this subsection are postponed to the next subsection.
\begin{lem}\label{lem6.4}

	\begin{enumerate}
		\item  There exists $c >0$ such that for any $n$ sufficiently great, $x \in \mathbb X$ and $a \in [n^{\frac{1}{2} -\varepsilon}, \theta_n n^{\frac{1}{2}}]   $,
			\begin{eqnarray} \label{eqn6.41}
		\left| \mathbb P_{x,a} (\tau >n)- \frac{2a}{ \sigma \sqrt{2 \pi n}} \right|   \le c  \frac{a\theta_n }{\sqrt{n}}.
		\end{eqnarray}
	
		\item  There exists $c >0$ such that for any  $x \in \mathbb X$, $a \ge n^{\frac{1}{2} -\varepsilon}$ and $n \ge 1$, 
		\begin{eqnarray} \label{eqn6.42}
		\mathbb P _{x,a} (\tau >n) \le c \frac{a}{\sqrt n}.
		\end{eqnarray}
	
	\end{enumerate}

	\end{lem}
Notice that for any $x \in \mathbb X, a \ge 0$ and $0 \le k \le n^{1-\varepsilon}$,
\begin{eqnarray}\label{eqn6.17}
P_n(x,a) \le P_{n-k}(x,a) \le P_{n -[n^{1-\varepsilon}]} (x,a).
\end{eqnarray}
By definition of $\nu_{n,\varepsilon}$ and (\ref{lem4.3}), as long as $A \le n^{\frac{1}{2}-\varepsilon}$, we have $\mathbb P _{x,a}$-a.s.
\vspace{-0.3cm}
\begin{eqnarray} \label{eqn6.16}
S_{\nu_{n,\varepsilon}} \ge M_{\nu_{n,\varepsilon}}- A  \ge 2 n^{\frac{1}{2}- \varepsilon} -A \ge n^{\frac{1}{2}- \varepsilon} .
\end{eqnarray}
Using (\ref{eqn6.41}) and (\ref{eqn6.17}), (\ref{eqn6.16}) with $\theta_n$ replaced by $\theta_n \left(\frac{n}{n - n^{1-\varepsilon}} \right) ^{\frac{1}{2}}$, for $n$ sufficiently great, on the event $ \left[ S_{\nu_{n,\varepsilon}} \le \theta_n n^{\frac{1}{2}}, \tau > \nu_{n,\varepsilon}, \nu_{n,\varepsilon} \le n^{1 - \varepsilon}\right]$,  we obtain $\mathbb P_{x,a}$-a.s.
$$
P_{n - \nu_{n,\varepsilon}} (X_{\nu_{n,\varepsilon}}, S_{\nu_{n,\varepsilon}}) = \frac{2(1+ o(1)) S_{\nu_{n,\varepsilon}}}{\sigma \sqrt {2 \pi n}}  .
$$
 Let 
 \begin{eqnarray}
 I_n' (x,a) &:=& \mathbb E_{x,a} \left[ S_{\nu_{n,\varepsilon}}; \tau > \nu_{n,\varepsilon}, \nu_{n,\varepsilon} \le n^{1-\varepsilon}  \right] \label{I'},   \\
 J_n' (x,a) &:=& \mathbb E_{x,a} \left[ S_{\nu_{n,\varepsilon}}; S_{\nu_{n,\varepsilon}} > \theta_n n^{\frac{1}{2}}, \tau > \nu_{n,\varepsilon}, \nu_{n,\varepsilon} \le n^{1-\varepsilon} \right]. \label{J'} 
 \end{eqnarray}
 Hence 
\begin{eqnarray*}
I_n(x,a) &=& \frac{2(1+o(1))}{\sigma \sqrt{2\pi n}} \mathbb E_{x,a} \left[ S_{\nu_{n,\varepsilon}}; S_{\nu_{n,\varepsilon}} \le \theta_n n^{\frac{1}{2}}, \tau > \nu_{n,\varepsilon}, \nu_{n,\varepsilon} \le n^{1-\varepsilon} \right]  \\
 &=& \frac{2(1+o(1))}{\sigma \sqrt{2\pi n}} \left[ I_n' (x,a) -J_n'(x,a) \right],\\
 J_n(x,a) &=& \frac{c(1+ o(1))}{\sqrt n} J_n'(x,a).
\end{eqnarray*}
Therefore (\ref{eqn6.13}) becomes 
$$
\Big| \mathbb P_{x,a} (\tau >n) - \frac{2(1+o(1))}{\sigma \sqrt{2\pi n}} I_n' (x,a) \Big|  \le  C \left( n^{-\frac{1}{2}} J_n'(x,a) \right) + C' \left( e^{-c_\varepsilon n^\varepsilon} \right).
$$
The first assertion of Theorem \ref{theo3} immediately follows by noticing that the term $J_n'$ is negligible and $\mathbb P_{x,a} (\tau >n)$ is dominated by the term $I_n'$ as shown in the lemma below.
\begin{lem} \label{lem6.5}
	
		$$\lim_{n \to +\infty} I_n'(x,a) = V(x,a) \quad \mbox{and} \quad \lim_{n \to +\infty} n^{2 \varepsilon} J_n' =0,
		$$
where $I_n'$ and $J_n'$ are defined in (\ref{I'}) and (\ref{J'}).
\end{lem}

{\bf (2)} By using Proposition \ref{theo2} (3), it suffices to prove $ \sqrt n \mathbb P_{x,a} (\tau >n) \le c(1+a)$ for $n$ great enough. For $n$ sufficiently great, using (\ref{eqn6.42}) and (\ref{eqn6.16}), we obtain  $\mathbb P_{x,a}$- a.s.
$$
P_{n -[n^{1-\varepsilon}]} (X_{\nu_{n,\varepsilon}}, S_{\nu_{n,\varepsilon}}) \le  c \frac{ S_{\nu_{n,\varepsilon}}}{\sqrt n}.
$$
Combined with (\ref{4.67}), it yields
\begin{eqnarray}\label{eqn6.43}
\mathbb P_{x,a} (\tau >n, \nu_{n,\varepsilon} \le n^{1-\varepsilon}) \le \frac{c}{\sqrt n} I_n'.
\end{eqnarray}
Since $\tau_a < T_{a+A}$ $\mathbb P$-a.s. and (\ref{5.15}), it follows that
$$
I'_n(x,a) \le \mathbb E_{x,a+A} [M_{\nu_{n,\varepsilon}}; T > \nu_{n,\varepsilon}, \nu_{n,\varepsilon} \le n^{1-\varepsilon} ] \le c(1+a+A).
$$
Hence (\ref{eqn6.43}) becomes
\begin{eqnarray} \label{eqn6.44}
\mathbb P_{x,a} (\tau >n, \nu_{n,\varepsilon} \le n^{1-\varepsilon}) \le \frac{c}{\sqrt n} (1+a+A).
\end{eqnarray}
Combining (\ref{eqn6.13}), (\ref{eqn6.14}) and (\ref{eqn6.44}), we obtain for $n$ great enough, 
$$
\mathbb P_{x,a} (\tau >n) \le e^{-c_\varepsilon n^\varepsilon} + \frac{c}{\sqrt n} (1+a+A) \le c'(1+a).
$$

{\rightline{$\square$}}

\subsection{Proof of Theorem \ref{theo4}}

Let us decompose $\mathbb P_{x,a} (S_n \le t \sqrt n|\tau >n)$ as follows:
\begin{eqnarray} \label{eqn7.11} 
\frac{\mathbb P_{x,a} (S_n \le t \sqrt n,\tau >n)} {\mathbb P_{x,a} (\tau >n)} = D_{n,1} + D_{n,2} + D_{n,3},
\end{eqnarray}
where 
\begin{eqnarray*}
	D_{n,1} &:=& \frac{ \mathbb P_{x,a} (S_n \le t \sqrt n,\tau >n, \nu_{n, \varepsilon} > n^{1-\varepsilon})}{\mathbb P_{x,a} (\tau >n)},
	\\
	D_{n,2} &:=& \frac{ \mathbb P_{x,a} (S_n \le t \sqrt n,\tau >n, S_n > \theta_n \sqrt n, \nu_{n, \varepsilon} \le n^{1-\varepsilon})}{\mathbb P_{x,a} (\tau >n)},
	\\
	D_{n,3} &:=& \frac{ \mathbb P_{x,a} (S_n \le t \sqrt n,\tau >n, S_n \le \theta_n \sqrt n, \nu_{n, \varepsilon} \le n^{1-\varepsilon})}{\mathbb P_{x,a} (\tau >n)}.
\end{eqnarray*}
Lemma \ref{lem5.4} and Theorem \ref{theo3} imply 
\begin{eqnarray}\label{eqn7.13}
\lim_{n \to +\infty} D_{n,1} =0.
\end{eqnarray}
 Theorem \ref{theo3} and Proposition \ref{theo2} (3) imply 
\begin{eqnarray*}
	D_{n,2}  &\le &  \frac{ \mathbb P_{x,a} ( \tau >n, S_n > \theta_n \sqrt n,  \nu_{n, \varepsilon} \le n^{1-\varepsilon})}{\mathbb P_{x,a} (\tau >n)} \notag \\
	&= & \frac{1}{\mathbb P_{x,a} (\tau >n)}  \mathbb E_{x,a} \Bigl[ P_{n- \nu_{n,\varepsilon}}(X_{\nu_{n,  \varepsilon}}, S_{\nu_{n,\varepsilon}}) ; \tau > \nu_{n,  \varepsilon}, S_{\nu_{n, \varepsilon}} > \theta_n \sqrt n, \nu_{n, \varepsilon} \le n^{1-\varepsilon} \Bigr] \notag  \\
	& \le& c \frac{ \mathbb E_{x,a} \Bigl[ 1 +S_{\nu_{n,\varepsilon}} ; \tau > \nu_{n,  \varepsilon}, S_{\nu_{n, \varepsilon}} > \theta_n \sqrt n, \nu_{n, \varepsilon} \le n^{1-\varepsilon} \Bigr] }{ \mathbb P_{x,a} (\tau >n)\sigma \sqrt{n - n^{1 - \varepsilon}} }  \notag \\
	&\le & c' \frac{\mathbb E_{x,a} \Bigl[  S_{\nu_{n,  \varepsilon} } ;  \tau > \nu_{n,  \varepsilon}, S_{\nu_{n, \varepsilon}} > \theta_n \sqrt n, \nu_{n, \varepsilon} \le n^{1-\varepsilon} \Bigr] + \mathbb P_{x,a}(\tau > \nu_{n,  \varepsilon})}{V(x,a)\sqrt{1 - n^{ - \varepsilon}} }.
\end{eqnarray*}
Since $\mathbb P_{x,a} (\tau <+ \infty) =1$ and $\mathbb P_{x,a} (\nu_{n, \varepsilon} <+ \infty) =0$,  Lemma \ref{lem6.5} yields
\begin{eqnarray} \label{eqn7.12}
\lim_{n \to +\infty} D_{n,2} =0.
\end{eqnarray}
Now we control $D_{n,3}$. Let $H_m (x,a) := \mathbb P_{x,a} ( S_m \le t \sqrt n, \tau >m)$.  We claim the following lemma and postpone its proof at the end of this section.
\begin{lem}\label{lem7.1}
	Let $\varepsilon \in (0,\varepsilon_0), t >0$ and $(\theta_n)_{n \ge 1}$ be a sequence such that $\theta_n \to 0$ and $\theta_n n^{\varepsilon/4} \to +\infty$ as $n \to +\infty$. Then for any $x \in \mathbb X$, $n^{1/2 - \varepsilon} \le a \le \theta_n \sqrt n$ and $1 \le k \le n^{1- \varepsilon}$, 
	\begin{eqnarray*}
		\mathbb P_{x,a} \Bigl(  S_{n-k} \le t \sqrt n, \tau >n-k  \Bigr) = \frac{2a}{\sigma^3 \sqrt{2 \pi n}}  \int_0^t u \exp \Bigl(-\frac{u^2}{2 \sigma^2} \Bigr)  du (1+o(1)).
	\end{eqnarray*}
\end{lem}
It is noticeable that on the event $ [\tau > k, S_k \le \theta_n \sqrt n, \nu_{n, \varepsilon} =k ] $, the random variable $H_{n-k} (X_k,S_k)$ satisfies the hypotheses of Lemma \ref{lem7.1}. Hence
\begin{eqnarray*}
	&& \mathbb P_{x,a} (S_n \le t \sqrt n,\tau >n, S_n \le \theta_n \sqrt n, \nu_{n, \varepsilon} \le n^{1-\varepsilon}) \\
	&=&  \mathbb E_{x,a} \Bigl[ H_{n- \nu_{n,  \varepsilon}} (X_{\nu_{n,  \varepsilon}}, S_{\nu_{n,  \varepsilon}} );  \tau > \nu_{n,  \varepsilon}, S_{\nu_{n,  \varepsilon}} \le \theta_n \sqrt n, \nu_{n, \varepsilon} \le n^{1-\varepsilon} \Bigr] \\
	&=& \sum_{k=1}^{[n^{1- \varepsilon}]} \mathbb E_{x,a} \Bigl[ H_{n- k} (X_k, S_k );  \tau > k, S_k \le \theta_n \sqrt n, \nu_{n, \varepsilon} =k \Bigr] \\
	&=& \frac{2 (1+o(1))}{\sigma^3 \sqrt{2 \pi n }}  \int_0^t u \exp \left( - u^2 \over 2 \sigma^2 \right) du \,  \mathbb E_{x,a} \Bigl[  S_{\nu_{n, \varepsilon}} ;  \tau > \nu_{n, \varepsilon}, S_{\nu_{n, \varepsilon}} \le \theta_n \sqrt n, \nu_{n, \varepsilon} \le n^{1- \varepsilon} \Bigr]  . 
\end{eqnarray*}
Lemma \ref{lem6.5}  yield as $n \to +\infty$, 
$$
\mathbb E_{x,a} \Bigl[  S_{\nu_{n, \varepsilon}} ;  \tau > \nu_{n, \varepsilon}, S_{\nu_{n, \varepsilon}} \le \theta_n \sqrt n, \nu_{n, \varepsilon} \le n^{1- \varepsilon} \Bigr] = V(x,a) (1+o(1)).
$$
Therefore, Theorem \ref{theo3} yields
\begin{eqnarray} \label{eqn7.14}
D_{n,3} &=& \frac{2 V(x,a) (1+o(1))}{ \mathbb P_{x,a} (\tau >n) \sigma^3 \sqrt{2 \pi n}} \int_0^t u \exp \left( - u^2 \over 2 \sigma^2 \right) du \notag \\
&=& \frac{1+o(1)}{\sigma^2} \int_0^t u \exp \left( - u^2 \over 2 \sigma^2 \right) du. 
\end{eqnarray}
The assertion of the theorem arrives by combining (\ref{eqn7.11}),  (\ref{eqn7.13}), (\ref{eqn7.12}) and (\ref{eqn7.14}).  

{\rightline{$\square$}} 

\subsection{Proof of Lemma {\ref{lem6.4}}}

{\bf (1)} Fix $\varepsilon >0$ and  let
$$
A_{n, \varepsilon}:=  \left[ \sup_{0 \le t \le 1} | S_{[nt]} - \sigma B_{nt}| \le  n^{\frac{1}{2} -2\varepsilon} \right].
$$
For any $x \in \mathbb X$, (\ref{eqn6.5}) implies $\mathbb P_x (A_{n, \varepsilon}^c) \le c_0 n^{-2\varepsilon}$. Denote $a^ {\pm} := a \pm n^{\frac{1}{2} -2\varepsilon} $ and notice that for $a \in [n^{\frac{1}{2} -\varepsilon}, \theta_n \sqrt n ]$,
\begin{eqnarray} \label{eqn6.9}
0 \le a^{\pm} \le 2 \theta_n \sqrt n.
\end{eqnarray}
Using (\ref{eqn6.4}) and (\ref{eqn6.9}), for any $x \in \mathbb X$ and $a \in [n^{\frac{1}{2} -\varepsilon}, \theta_n \sqrt n]$, we obtain

\begin{eqnarray}\label{eqn6.10}
 - \frac{c a^{\pm} \theta_n}{\sqrt n} \pm  \frac{2n^{-2 \varepsilon}}{\sigma \sqrt {2 \pi}} \le \mathbb P_x (\tau_{a^\pm}^{bm} >n) -  \frac{2a}{\sigma \sqrt{2 \pi n}} \le   \frac{c a^{\pm} \theta_n}{\sqrt n} \pm  \frac{2n^{-2 \varepsilon}}{\sigma \sqrt {2 \pi}}.
\end{eqnarray}

For any $a \ge n^{\frac{1}{2} -\varepsilon}$, we have 
$
\left[\tau_{a^-}^{bm} >n \right] \cap A_{n, \varepsilon}^c \subset \left[\tau_a >n \right] \cap A_{n, \varepsilon}^c \subset \left[ \tau_{a^+}^{bm} >n \right] \cap A_{n, \varepsilon}^c,
$
which yields 
$$
\mathbb P_x (\tau_{a^-}^{bm} >n) - \mathbb P_x (A^c_{n, \varepsilon}) \le \mathbb P_x (\tau_a >n) \le \mathbb P_x (\tau_{a^+}^{bm} >n) + \mathbb P_x (A^c_{n, \varepsilon})
$$
 for any $x \in \mathbb X$. It follows that 
\begin{eqnarray}\label{eqn6.8}
\left\{ \begin{array}{l}
\mathbb P_x (\tau_a >n) - \mathbb P_x (\tau_{a^+}^{bm} >n)  \le c_0 n^{-2\varepsilon}, \\
\mathbb P_x (\tau_{a^-}^{bm} >n) -  \mathbb P_x (\tau_a >n)   \le c_0 n^{-2\varepsilon}.
\end{array} \right.
\end{eqnarray}
The fact that $\theta_n n^{\varepsilon/4} \to +\infty$ yields for $n$ great enough
\begin{eqnarray}\label{eqn6.11}
\theta_n \frac{a}{\sqrt n} \ge \frac{n^{\frac{1}{2} -\varepsilon}}{n^\varepsilon \sqrt n} = n^{-2\varepsilon}.
\end{eqnarray}
From (\ref{eqn6.10}), (\ref{eqn6.8}) and (\ref{eqn6.11}), it follows that for any  $a \in [n^{\frac{1}{2} -\varepsilon}, \theta_n \sqrt n]$, 

$$
\left| \mathbb P_x (\tau_a > n) - \frac{2a}{\sigma \sqrt{2 \pi n}}   \right| \le c(1+ \theta_n) n^{-2\varepsilon} + c_1\frac{\theta_n a}{\sqrt n} \le c_2 \frac{\theta_n a}{\sqrt n}.
$$

{\bf (2)} For $n$ great enough, condition $a \ge n^{\frac{1}{2} -\varepsilon}$ implies $a^+ \le 2a$. From  (\ref{eqn6.3}) and (\ref{eqn6.8}) , since $n^{-2\varepsilon} \le \frac{a}{\sqrt n}$, for any $x \in \mathbb X$,
$$
\mathbb P_x (\tau_a >n) \le c \frac{a}{\sigma \sqrt n} + c_0 n^{-2\varepsilon} \le c_1 \frac{a}{\sqrt n}. 
$$

{\rightline{$\square$}}

\subsection{Proof of Lemma \ref{lem6.5}}

{\bf (1)} We prove that $\displaystyle \lim_{n \to +\infty} \mathbb E_{x,a} \left[ M_{\nu_{n,\varepsilon}}; \tau > \nu_{n,\varepsilon}, \nu_{n,\varepsilon} \le n^{1-\varepsilon}  \right] = V(x,a)$. Then, the assertion arrives by using (\ref{lem4.3}) and taking into account that $\mathbb P_x (\tau_a < +\infty) =1$ and  $\displaystyle \mathbb P_x (\lim_{n \to +\infty} \nu_{n,\varepsilon} = +\infty) =1$. 
For $x \in \mathbb X$  and $a  \ge 0$, we obtain
\begin{eqnarray}
\mathbb E_{x,a} \left[ M_{\nu_{n,\varepsilon}}; \tau > \nu_{n,\varepsilon}, \nu_{n,\varepsilon} \le n^{1-\varepsilon}  \right]
&=& \mathbb E_{x,a} \left[ M_{\nu_{n,\varepsilon} \wedge [n^{1 -\varepsilon}]}; \tau > \nu_{n,\varepsilon}, \nu_{n,\varepsilon} \le n^{1-\varepsilon}  \right] \notag \\
 &=& \mathbb E_{x,a} \left[ M_{\nu_{n,\varepsilon} \wedge [n^{1 -\varepsilon}]}; \tau > \nu_{n,\varepsilon} \wedge [n^{1 -\varepsilon}]  \right]  \notag \\
 & &  \,\, - \mathbb E_{x,a} \left[ M_{\nu_{n,\varepsilon} \wedge [n^{1 -\varepsilon}]}; \tau > \nu_{n,\varepsilon} \wedge [n^{1 -\varepsilon}], \nu_{n,\varepsilon} > n^{1-\varepsilon}  \right] \notag.
\end{eqnarray}
 By using Lemma \ref{lem5.5},
 \begin{eqnarray*} \label{eqm6.25}
 \mathbb E_{x,a} \left[ M_{\nu_{n,\varepsilon} \wedge [n^{1 -\varepsilon}]}; \tau > \nu_{n,\varepsilon} \wedge [n^{1 -\varepsilon}], \nu_{n,\varepsilon} > n^{1-\varepsilon}  \right] \le  c(1+a) e^{-c_\varepsilon n^\varepsilon}.
 \end{eqnarray*}
Using the facts that $(M_n)_{n \ge 0}$ is a martingale and $\displaystyle \mathbb P_x \left(\lim_{n \to +\infty}\nu_{n,\varepsilon} = + \infty\right) =1$, we obtain 
\begin{eqnarray*}
  \lim_{n \to +\infty} \mathbb E_{x,a} \left[ M_{\nu_{n,\varepsilon}}; \tau > \nu_{n,\varepsilon}, \nu_{n,\varepsilon} \le n^{1-\varepsilon}  \right]  &=& \lim_{n \to +\infty}  \mathbb E_{x,a} \left[ M_{\nu_{n,\varepsilon} \wedge [n^{1 -\varepsilon}]}; \tau > \nu_{n,\varepsilon} \wedge [n^{1 -\varepsilon}]  \right]  \notag \\
 &=& a- \lim_{n \to +\infty} \mathbb E_{x,a} \left[ M_{\nu_{n,\varepsilon} \wedge [n^{1 -\varepsilon}]}; \tau \le \nu_{n,\varepsilon} \wedge [n^{1 -\varepsilon}]  \right]   \notag \\
 &=& a-  \lim_{n \to +\infty} \mathbb E_{x,a} \left[ M_\tau; \tau \le \nu_{n,\varepsilon} \wedge [n^{1 -\varepsilon}]  \right]  \notag \\
 &=& a - \mathbb E _{x,a} [M_\tau]= V(x,a).
\end{eqnarray*}

{\bf (2)} Let $b = a +A$. Remind that $M^*_n =  \displaystyle \max_{1 \leq k \leq n} |M_k|$. We obtain
\begin{eqnarray*}
	\mathbb E_{x,a} \left[ S_{\nu_{n,\varepsilon}}; S_{\nu_{n,\varepsilon}} > \theta_n n^{\frac{1}{2}}, \tau > \nu_{n,\varepsilon}, \nu_{n,\varepsilon} \le n^{1-\varepsilon} \right] &\le & \mathbb E_{x,b} \left[ M_{\nu_{n,\varepsilon}}; M_{\nu_{n,\varepsilon}} > \theta_n n^{\frac{1}{2}}, \nu_{n,\varepsilon} \le n^{1-\varepsilon} \right] \\
	&\le & \mathbb E_{x,b} \left[ M^*_{[n^{1-\varepsilon}]};  M^*_{[n^{1-\varepsilon}]} > \theta_n n^{\frac{1}{2}} \right] .
\end{eqnarray*}
Since $\theta_n n^{\varepsilon/4} \to +\infty$ as $n \to +\infty$, it suffices to prove that for any $\delta >0$, $ x \in \mathbb X$ and $b \in \mathbb R$,
\begin{eqnarray*}\label{eqn6.32}
	\lim_{n \to +\infty} n^{2 \varepsilon} \mathbb E _x \left[b+M^*_n; M^*_n > n^{\frac{1}{2} +\delta} \right] =0.
\end{eqnarray*}
Obviously, by (\ref{eqn5.8}),
\begin{eqnarray*}
	\mathbb E _x \left[b+M^*_n; M^*_n > n^{\frac{1}{2} +\delta} \right] &\le& b \mathbb P _x \left( M^*_n > n^{\frac{1}{2} +\delta} \right) +\mathbb E_x \left[M^*_n; M^*_n > n^{\frac{1}{2} +\delta}  \right] \\
	&=&  \left( b+ n^{\frac{1}{2} +\delta} \right) \mathbb P _x \left( M^*_n > n^{\frac{1}{2} +\delta} \right) + \int^{+\infty}_{n^{\frac{1}{2} +\delta}} \mathbb P_x (M^*_n >t) dt\\
	&\le& c\left( b+n^{\frac{1}{2}+\delta}\right) n^{-p\delta} + c n ^{-p\delta + \frac{1}{2} +\delta}.
\end{eqnarray*}
Since $p$ can be taken arbitrarily great, it follows that $\displaystyle \lim_{n \to +\infty} n^{2 \varepsilon} J_n' =0$.

{\rightline{$\square$}}

\subsection{Proof of lemma \ref{lem7.1}}

Recall that $a^\pm = a \pm n^{1/2 - 2 \varepsilon}$ and denote $t^\pm = t \pm 2 n^{-2 \varepsilon}$. For any $1 \le k \le n^{1- \varepsilon}$, 
$$
\{\tau_{a^-}^{bm} \} \cap A_{n, \varepsilon} \subset \{ \tau_a > n-k\} \cap A_{n, \varepsilon} \subset \{\tau_{a^+}^{bm} \} \cap A_{n, \varepsilon}
$$
and 
$$
\{a^-+\sigma B_{n-k} \le t^-\sqrt n \} \cap A_{n, \varepsilon} \subset \{a+S_{n-k} \le t\sqrt n \} \cap A_{n, \varepsilon} \subset \{a^+ +\sigma  B_{n-k} \le t^+\sqrt n \} \cap A_{n, \varepsilon},
$$
which imply
\begin{eqnarray}\label{eqn7.2}
&& \mathbb  P_x (\tau_{a^-}^{bm} >n-k, a^- +\sigma  B_{n-k} \le t^-\sqrt n) - \mathbb P_x (A^c_{n, \varepsilon}) \notag \\
&&\qquad \qquad \qquad \le  \mathbb P_x (\tau_a >n-k, a+ S_{n-k} \le t \sqrt n) \le  \\
 &&   \qquad \qquad \qquad \qquad \qquad  \mathbb P_x (\tau_{a^+}^{bm} >n-k, a^+ +\sigma  B_{n-k} \le t^+\sqrt n) +\mathbb P_x (A^c_{n, \varepsilon}).\notag 
\end{eqnarray}
Moreover, by Lemma \ref{lem6.1}, we obtain
\begin{eqnarray}\label{eqn7.8}
\mathbb P_x  \Bigl(\tau_{a^+}^{bm} >n-k, a^+ +\sigma  B_{n-k} \le t^+\sqrt n  \Bigr) = \frac{2a}{\sigma^3 \sqrt{2 \pi n}} \int_0^t u \exp \Bigl(-\frac{u^2}{2 \sigma^2} \Bigr)  du (1+o(1))
\end{eqnarray}
and similarly,
\begin{eqnarray}\label{eqn7.9}
\mathbb P_x  \Bigl(\tau_{a^-}^{bm} >n-k, a^- +\sigma  B_{n-k} \le t^-\sqrt n  \Bigr) = \frac{2a}{\sigma^3 \sqrt{2 \pi n}} \int_0^t u \exp \Bigl(-\frac{u^2}{2 \sigma^2} \Bigr)  du (1+o(1)).
\end{eqnarray}
Therefore, from  (\ref{eqn7.2}), (\ref{eqn7.8}), (\ref{eqn7.9}) and $\mathbb P_x (A^c_{n, \varepsilon})\le c n^{-2\varepsilon}$, it follows that

\begin{eqnarray*}
	\mathbb P_x  \Bigl(\tau_a >n-k, a + S_{n-k} \le t \sqrt n  \Bigr) = \frac{2a}{\sigma^3 \sqrt{2 \pi n}}  \int_0^t u \exp \Bigl(-\frac{u^2}{2 \sigma^2} \Bigr)  du (1+o(1)).
\end{eqnarray*}

{\rightline{$\square$}}

\section{On conditions C1-C3 of Theorem 2.1 in \cite{GLP1}}

Let $k_{gap}, M_1, M_2 \in \mathbb N$ and $j_0 < \ldots < j_{M_1 +M_2}$ be natural numbers. Denote $a_{k + J_m} = \sum_{l \in  J_m} a_{k+l}$, where $J_m = [j_{m-1}, j_m), m= 1, \ldots, M_1 +M_2$ and $k \ge 0$. Consider the vectors $\bar {a}_1 = (a_{J_1}, \ldots, a_{J_{M_1}})$ and $\bar {a}_2 = (a_{k_{gap}+J_{M_1 +1}}, \ldots, a_{k_{gap}+J_{M_1 + M_2}})$. Denote by $\phi_x (s,t) = \mathbb E e^{is\bar {a}_1 + it \bar {a}_2 }$, $\phi_{x,1} (s) = \mathbb E e^{is\bar {a}_1}$ and $\phi_{x,2} (s) = \mathbb E e^{it\bar {a}_2}$ the characteristic functions of $(\bar {a}_1, \bar {a}_2 )$, $\bar {a}_1$ and $\bar {a}_2 $, respectively. For the sake of brevity, we denote $\phi_1(s) = \phi_{x,1}(s), \phi_2(t) = \phi_{x,2}(t)$ and $\phi(s,t) = \phi_{x}(s,t)$.

 We first check that conditions C1-C3 hold and then prove the needed lemmas.

\subsection{Statement and proofs of conditions C1-C3}

{\bf  C1}: There exist positive constants $\varepsilon_0 \le 1, \lambda_0, \lambda_1, \lambda_2$ such that for any $k_{gap} \in \R, M_1, M_2 \in \N$, any sequence $j_0 < \ldots < j_{M_1 + M_2}$ and any $s \in \R ^{M_1}, t \in \R ^{M_2}$ satisfying $|(s,t)|_\infty \le \varepsilon_0$, 
$$
|\phi(s,t) - \phi_1(s) \phi_2(t)| \le \lambda_0 \exp (-\lambda_1 k_{gap}) \left( 1+ \max _{m= 1, \ldots, M_1+M_2} card(J_m)\right)^{\lambda_2(M_1+M_2)}.
$$
{\bf  C2}: There exists a positive constant $\delta $  such that 
$\sup_{n \ge 0} | a_n|_{L^{2 + 2 \delta}}  < +\infty$.\\
{\bf  C3}: There exist a positive constant $C$ and a positive number $\sigma$ such that for any $\gamma >0$, any $x \in \X$ and any $n \ge 1$, 
$$
\sup_{m \ge 0} \left| n^{-1} Var_{\P_x} \left( \sum_{i=m}^{m+n-1} a_i\right) - \sigma^2 \right| \le C n^{-1 + \gamma}.
$$

\begin{prop} \label{C1}
	Condition 1 is satisfied under hypotheses P1-P5.
\end{prop}
{\bf Proof.}
First, we prove the following lemma.
\begin{lem} \label{lem3} There exist two positive constants $C$ and $\kappa$ such that $0< \kappa <1$ and 
	$$|\phi(s,t) - \phi_1(s) \phi_2(t)| \le C C_P^{M_1 +M_2} \kappa^{k_{gap}}, $$
	where $C_P$ is defined in Proposition \ref{spectre}.
\end{lem}
{\bf Proof.}
In fact, the characteristic functions of the ramdom variables $\bar{a}_1, \bar{a}_2$ and $(\bar{a}_1,\bar{a}_2)$ can be written in terms of operator respectively as follows:
\begin{eqnarray}
	\phi_1(s) &=& \E _x [e^{is \bar a_1}] = P ^{j_0 -1} P_{s_1}^{|J_1|} \ldots P_{s_{M_1}}^{|J_{M_1}|} {\bf 1} (x) ,\notag \\
	\phi_2(t) &=& \E _x [e^{it \bar a_2}] = P ^{k_{gap} +j_{M_1} -1} P_{t_1}^{|J_{M_1+1}|} \ldots P_{t_{M_2}}^{|J_{M_1 +M_2}|} {\bf 1} (x),  \\
	\phi(s,t) &=& \E _x [e^{is \bar a_1+ it \bar a_2}] = P ^{j_0 -1} P_{s_1}^{|J_1|} \ldots P_{s_{M_1}}^{|J_{M_1}|} P ^{k_{gap}} P_{t_1}^{|J_{M_1+1}|} \ldots P_{t_{M_2}}^{|J_{M_1 +M_2}|} {\bf 1} (x). \notag
\end{eqnarray}

Now we decompose $\phi(s,t)$ into the sum of $\phi_\Pi (s,t)$ and $\phi_R (s,t)$ by using the spectral decomposition $P = \Pi + R$ in Proposition \ref{spectre}, where 
\begin{eqnarray}
	\phi_\Pi (s,t) &=&  P ^{j_0 -1} P_{s_1}^{|J_1|} \ldots P_{s_{M_1}}^{|J_{M_1}|} \Pi P_{t_1}^{|J_{M_1+1}|} \ldots P_{t_{M_2}}^{|J_{M_1 +M_2}|} {\bf 1} (x), \notag \\
	\phi_R (s,t)  &=&  P ^{j_0 -1} P_{s_1}^{|J_1|} \ldots P_{s_{M_1}}^{|J_{M_1}|} R ^{k_{gap}} P_{t_1}^{|J_{M_1+1}|} \ldots P_{t_{M_2}}^{|J_{M_1 +M_2}|} {\bf 1} (x). \notag 
\end{eqnarray}
Since $\Pi(\varphi) = \nu (\varphi) {\bf 1}$ for any $\varphi \in L$ and $P_t$ acts on $L$, we obtain 
$$
\phi_\Pi (s,t) =  P ^{j_0 -1} P_{s_1}^{|J_1|} \ldots P_{s_{M_1}}^{|J_{M_1}|}  {\bf 1} (x) \nu\left(P_{t_1}^{|J_{M_1+1}|} \ldots P_{t_{M_2}}^{|J_{M_1 +M_2}|} {\bf 1}\right).
$$
Then setting $\psi_2(t) = \nu(P_{t_1}^{|J_{M_1+1}|} \ldots P_{t_{M_2}}^{|J_{M_1 +M_2}|} {\bf 1}) $ yields
\begin{align}
	\phi(s,t) &=  \phi_1(s) \psi_2(t) + \phi_R(s,t)  \notag \\
	&= \phi_1(s) \phi_2(t) + \phi_1(s) [ \psi_2(t) - \phi_2(t)] + \phi_R(s,t), \notag 
\end{align}
which implies 
\begin{eqnarray}\label{46}
	|\phi (s,t)- \phi_1(s) \phi_2(t) | \le |\phi_1(s)  | |\psi_2(t) - \phi_2(t) | + |\phi_R(s,t) |.
\end{eqnarray}
On the one hand, we can see that  $|\phi_1(s)| = \left| \left(  P ^{j_0 -1} P_{s_1}^{|J_1|} \ldots P_{s_{M_1}}^{|J_{M_1}|} {\bf 1} \right) (x) \right| \le C_P^{1+M_1}$ and $|\phi_R(s,t)|  \le C_P^{1+M_1 + M_2} C_R \kappa^{k_{gap}}$.
On the other hand, since $\nu $ is $P$-invariant measure and $(\nu - \delta_x) ({\bf 1} ) =0$, by using again the expression $P = \Pi +R$, we obtain
\begin{eqnarray}
	|\psi_2(t) - \phi_2(t) | &=& \left| (\nu - \delta_x) \left(  P ^{k_{gap} +j_{M_1} -1} P_{t_1}^{|J_{M_1+1}|} \ldots P_{t_{M_2}}^{|J_{M_1 +M_2}|} {\bf 1} \right) \right| \notag \\
	&\le& \left| (\nu - \delta_x) \left( \Pi P_{t_1}^{|J_{M_1+1}|} \ldots P_{t_{M_2}}^{|J_{M_1 +M_2}|} {\bf 1} \right) \right| \notag \\
	&&  \qquad \qquad  \qquad  + \left|  (\nu - \delta_x) \left(  R ^{k_{gap} +j_{M_1} -1} P_{t_1}^{|J_{M_1+1}|} \ldots P_{t_{M_2}}^{|J_{M_1 +M_2}|} {\bf 1} \right)\right| \notag \\
	&=&  \left| (\nu - \delta_x) ({\bf 1} ) \nu \left( P_{t_1}^{|J_{M_1+1}|} \ldots P_{t_{M_2}}^{|J_{M_1 +M_2}|} {\bf 1} \right) \right|  \notag \\
	&&   \qquad \qquad \qquad + \left|  (\nu - \delta_x) \left(  R ^{k_{gap} +j_{M_1} -1} P_{t_1}^{|J_{M_1+1}|} \ldots P_{t_{M_2}}^{|J_{M_1 +M_2}|} {\bf 1} \right)\right| \notag \\
	&=& \left| (\nu - \delta_x) \left(  R ^{k_{gap} +j_{M_1} -1} P_{t_1}^{|J_{M_1+1}|} \ldots P_{t_{M_2}}^{|J_{M_1 +M_2}|} {\bf 1} \right) \right| \notag \\
	&\le& C C_P^{M_2}  \kappa ^{k_{gap} +j_{M_1} -1} .
\end{eqnarray}
Therefore, (\ref{46}) follows.

\rightline{$\square$}
Second, let $\lambda_2 = \max \{1 , \log_2 C_P \}$. Since $\displaystyle \max_{m=1, \ldots, M_1 +M_2} card (J_m) \ge 1$, we obtain 
$$
C_P^{M_1 +M_2} \le 2^{\lambda_2 ( M_1 +M_2)} \le \left( 1+ \max_{m=1, \ldots, M_1 +M_2} card (J_m)\right) ^{\lambda_2 ( M_1 +M_2)},
$$
which implies that 
$$
|\phi(s,t) - \phi_1(s) \phi_2(t)| \le C \kappa^{k_{gap}}  \left( 1+ \max_{m=1, \ldots, M_1 +M_2} card (J_m)\right) ^{\lambda_2 ( M_1 +M_2)}.
$$

Finally, let $\lambda_0 =C$ and  $\lambda_1 = - \log \kappa$. Then the assertion arrives.

\rightline{$\square$}

\begin{prop} \label{C2}
	Condition 2 is satisfied under hypotheses P1-P5.
\end{prop}
{\bf Proof.} Condition P1 implies that there exists $\delta_0>0$ such that $\mathbb E [N(g)^{\delta_0}] < + \infty $ and since $\displaystyle \mathbb E [N(g)^{\delta_0}] = \mathbb E [\exp (\delta_0 \log N(g))] = \sum_{k=0}^{+\infty} \frac{\delta_0^k}{k!} \mathbb E [(\log N(g))^k] $, we obtain 
$\mathbb E |a_n|^k \le \mathbb E [(\log N(g))^k] < +\infty $
for any $n \ge 0$ and any $k \ge 0$.

\rightline{$\square$}

\begin{prop} \label{C3}
	Condition 3 is satisfied under hypotheses P1-P5. More precisely, there exists a positive  constant $\sigma $  such that for any $x \in \mathbb X$ and any $n \ge 1$,
	\begin{eqnarray} \label{C3}
		\sup_{m \ge 0} \left| Var_{\mathbb P_x} \left( \sum_{k=m}^{m+n-1} a_k \right) - n \sigma^2 \right| < +\infty.
	\end{eqnarray}
\end{prop}
{\bf Proof.} For any integer $m,n \ge 0$, we denote $S_{m,n} = \sum_{k=m}^{m+n-1}a_k$, $V_x (X) = Var_{\mathbb P_x} (X) = \mathbb E_x (X^2)- (\mathbb E_x X)^2$ and $Cov_x(X,Y) = Cov_{\mathbb P_x} (X,Y)$. Then 
\begin{eqnarray} \label{e1}
	V_x (S_{m,n}) &=& \sum_{k=m}^{m+n-1} V_x (a_k ) + 2 \sum_{k=m}^{m+n-1} \sum_{l=1}^{m+n-k-1} Cov_x (a_k, a_{k+l})
\end{eqnarray}
and (\ref{C3}) becomes $\sup_{m \ge 0} |V_x (S_{m,n}) - n \sigma^2 | < +\infty$. We claim two lemmas and postpone their proofs until the end of this section.
\begin{lem} \label{lem1}
	There exist $C>0$ and $0 < \kappa <1$ such that for any $x \in \mathbb X$, any $k \ge 0$ and any $l \ge 0$,
	\begin{eqnarray} \label{e3}
		\left| Cov_x (a_k, a_{k+l})\right| \le C \kappa^l.
	\end{eqnarray}
\end{lem}

\begin{lem} \label{lem2}
	There exist $C>0$, $0 < \kappa <1$ and a sequence $(s_n)_{n \ge 0}$ of real numbers  such that for any $x \in \mathbb X$, any $k \ge 0$ and any $l \ge 0$,
	\begin{align} 
		\left| Cov_x (a_k, a_{k+l}) - s_l \right| & \le C \kappa^k, \label{e4} \\
		|s_l| \le C \kappa^l. \label{e5}
	\end{align}

\end{lem}
For the first term of the right side of (\ref{e1}), by combining Lemma \ref{lem1} and Lemma \ref{lem2}, we obtain 
\begin{eqnarray}\label{e6}
	\left| Cov_x (a_k, a_{k+l}) - s_l \right| \le C \kappa^{\max\{k,l \}} .
\end{eqnarray}
Inequality (\ref{e4}) implies $ |V_x (a_k) -s_0 | \le C \kappa^k$, which yields for any integer $m,n \ge 0$,
\begin{eqnarray}\label{e7}
	\left|\sum_{k=m}^{m+n-1} V_x(a_k) - n s_0 \right| \le \sum_{k=m}^{m+n-1} |V_x(a_k)- s_0 | \le C \sum_{k=m}^{m+n-1} \kappa^k \le {C \over 1-\kappa} < +\infty.
\end{eqnarray}
For the second term of the right side of (\ref{e1}), we can see that
\begin{eqnarray}\label{e8}
	& &\left| \sum_{k=m}^{m+n-1} \sum_{l=1}^{m+n-k-1} Cov_x(a_k, a_{k+l}) - \sum_{k=m}^{m+n-1} \sum_{l=1}^{+\infty} s_l \right| \notag \\
	&\le& \sum_{k=m}^{m+n-1} \sum_{l=1}^{m+n-k-1} \left| Cov_x(a_k, a_{k+l}) - s_l \right| + \sum_{k=m}^{m+n-1} \sum_{l=m+n-k}^{+\infty} |s_l| \notag\\
	&=& \Sigma_1(x, m ,n) + \Sigma_2(x, m ,n).
\end{eqnarray}
On the one hand, by (\ref{e4}) and (\ref{e6}), we can see that for any $x \in \X$, any $m\ge 0$ and any $n \ge 1$,
\begin{eqnarray}\label{e9}
	\Sigma_1(x, m ,n) &\le& \sum_{k=0}^{+\infty} \sum_{l=1}^{k}C \kappa^k + \sum_{k=0}^{+\infty} \sum_{l=k+1}^{+\infty} C \kappa ^l \notag \\
	&\le& \sum_{k=0}^{+\infty} C k \kappa^k + \sum_{k=0}^{+\infty} C  {\kappa^{k+1} \over 1 - \kappa} < +\infty.
\end{eqnarray}
Similarly, on the other hand, by (\ref{e5}) we obtain for any $x \in \X$, any  $m\ge 0$ and any $n \ge 1$,
\begin{eqnarray}\label{e10}
	\Sigma_2(x, m ,n)  \le  \sum_{k=0}^{n-1} \sum_{l=n-k}^{+\infty} C \kappa^{l} \le  {C \over (1 - \kappa)^2} < +\infty.
\end{eqnarray}
Combining (\ref{e1}),(\ref{e7}),(\ref{e8}),(\ref{e9}) and (\ref{e10}), we obtain 
\begin{eqnarray}
	\sup_{m \ge 0} \left|V_x(S_{m,n}) - n \sum_{l=0}^{+\infty} s_l \right| < +\infty.
\end{eqnarray}

In fact, by using Lemma 2.1 in \cite{LPP}, Theorem 5 in \cite{H2} implies that the sequence $(\frac{S_n}{\sqrt n} )_{n \ge 1}$ converges weakly to a normal law with variance $\sigma^2$. Meanwhile,  under hypothesis P2, Corollary 3 in \cite{H2} implies that the sequence $(|R_n|)_{n \ge 1}$ is not tight and thus $\sigma^2 >0$, see \cite{H2} for the definition and basic properties.    Therefore, we can see that $Var_x S_n \sim n \sigma^2$ with $\sigma^2 >0$, which yields $\sum_{l=0}^{+\infty} s_l = \sigma^2$.

\rightline{$\square$}

\subsection{Proof of Lemma \ref{lem1}}

Let $g(x) = \left\{ \begin{array}{l}
x \ \ \mbox{if} \  |x| \le 1,\\
0\ \ \mbox{if} \ |x| > 2.
\end{array} \right.$ such that $g$ is $C^\infty$ on $\R$ and $|g(x)| \le |x|$ for any $x \in \R$. Then $g \in L^1(\R)\cap C^1_{c}(\R)$. Therefore, the Fourier transform of $g$ is $\hat g$ defined as follows: 
$$
\hat g (t) := \int_{\R}e^{-itx} g(x) dx,
$$
and the Inverse Fourier Theorem yields
$$
g(x) = \frac{1}{2 \pi}\int_{\R} e^{itx} \hat g(t) dt.
$$
Let $g_T(x) := Tg({x \over T})$ for any $T >0$. Then $|\hat g_T|_1 = T |\hat g|_1 < +\infty$.  Let $h_T (x,y) = g_T (x) g_T(y)$. Then $\hat h_T(x,y) = \hat g_T(x) \hat g_T(y)$. Let $V$ and $V'$ be two i.i.d. random variables with mean $0$, independent of $a_l$ for any $l \ge 0$ whose characteristic functions have the support included in the interval $[-\varepsilon_0, \varepsilon_0]$ for  $\varepsilon_0 $ defined in C1. Assume that  $\E |V|^n < +\infty$ for any $n >0$. Let $Z_k = a_k +V$ and $Z'_{k+l} = a_{k+l} +V'$ and denote by  $\widetilde \phi _1(s), \widetilde \phi_2 (t)$ and $\widetilde \phi (s,t)$ the characteristic functions of $Z_k, Z'_{k+l}$ and $(Z_k, Z'_{k+l})$, respectively. 

We use the same notations introduced at the beginning of this section by setting $\phi_1(s) = \E _x [e^{isa_k}], \phi_2(t) = \E _x [e^{ita_{k+l}}]$ and $\phi(s,t) = \E _x [e^{is a_k + it a_{k+l}}]$. We also denote $\varphi$ the characteristic function of $V$, that yields
\begin{eqnarray}\label{1.23}
	\widetilde \phi_1(s) &=& \E [e^{isZ_k}] = \E [e^{isa_k}] \E [e^{isV}] =\phi_1(s) \varphi(s) ,\notag \\
	\widetilde \phi_2(t) &=& \E [e^{itZ'_{k+l}}] = \E [e^{ita_{k+l}}] \E [e^{itV'}] =\phi_2(t) \varphi(t),  \\
	\widetilde \phi(s,t) &=& \E [e^{isZ_k + itZ'_{k+l} }] = \E [e^{isa_k + ita_{k+l} }] \E [e^{isV}] \E [e^{itV'}] =\phi(s,t) \varphi(s) \varphi(t) .\notag
\end{eqnarray}
Then we can see that $\widetilde \phi_1$ and $\widetilde \phi_2$  have the support in $[-\varepsilon_0, \varepsilon_0]$. We perturb $a_k$ and $a_{k+l}$ by adding the random variables $V$ and $V'$  with mean $0$ and the support of their characteristic functions are on $[-\varepsilon_0, \varepsilon_0]$.
We explicit the quantity $Cov_x (a_k, a_{k+l})$:
\begin{eqnarray} \label{11}
	Cov_x (a_k, a_{k+l}) = \E _x [a_k, a_{k+l}] - \E _x a_k  \E _x a_{k+l}.
\end{eqnarray}
On the one hand, we can see that 
\begin{eqnarray} \label{12}
	\E _x [a_k a_{k+l}] = \E _x [Z_k Z'_{k+l}] &=& \E _x [h_T(Z_k; Z'_{k+l})] +\E _x [Z_k Z'_{k+l}] - \E _x [h_T(Z_k; Z'_{k+l})] \notag \\
	&=&\frac{1}{(2 \pi)^2} \E _x \int \int  \hat h _T (s,t) e^{isZ_k + it Z'_{k+l}} ds dt + R_0 \notag \\
	&=&\frac{1}{(2 \pi)^2}  \int \int \hat h _T (s,t)  \E _x \left[ e^{isZ_k + it Z'_{k+l}} \right] ds dt + R_0 \notag \\
	&=&\frac{1}{(2 \pi)^2}  \int \int \hat h _T (s,t)  \widetilde \phi (s,t) ds dt + R_0, \notag \\
\end{eqnarray}
where $R_0 = \E _x [Z_k Z'_{k+l}] - \E _x [h_T(Z_k; Z'_{k+l})]$. On the other hand, we obtain
\begin{eqnarray} \label{13}
	\E _x a_k = \E _x Z_k  &=& \E _x g_T(Z_k)  +\E _x Z_k  - \E _x g_T(Z_k) \notag \\
	&=& \frac{1}{2 \pi}  \int \hat g _T (s) \widetilde \phi_1 (s)ds + R_1, 
\end{eqnarray}
where $R_1 =\E _x Z_k  - \E _x g_T(Z_k) $ and
\begin{eqnarray} \label{14}
	\E _x a_{k+l} = \E _x Z'_{k+l}  &=& \E _x g_T(Z'_{k+l})  +\E _x Z'_{k+l}  - \E _x g_T(Z'_{k+l}) \notag \\
	&=& \frac{1}{2 \pi}  \int \hat g _T (t) \widetilde \phi_2 (t)dt + R_2, 
\end{eqnarray}
where $R_2 = \E _x Z'_{k+l}  - \E _x g_T(Z'_{k+l}) $.
From (\ref{11}), (\ref{12}), (\ref{13}) and (\ref{14}), since $\hat h_T(s,t) = \hat g_T(s) \hat g_T(t)$, we obtain 
\begin{eqnarray}\label{1.28}
	Cov_x (a_k, a_{k+l}) &=& \frac{1}{(2 \pi)^2}  \int \int \hat h _T (s,t)  \widetilde \phi (s,t) ds dt + R_0   \notag \\
	& & \qquad  -\left(\frac{1}{2 \pi}  \int \hat g _T (s) \widetilde \phi_1 (s)ds + R_1  \right) \left( \frac{1}{2 \pi}  \int \hat g _T (t) \widetilde \phi_2 (t)dt + R_2 \right) \notag \\
	&=& \frac{1}{(2 \pi)^2}  \int \int \hat h _T (s,t)   \left[ \widetilde \phi (s,t) - \widetilde \phi_1(s) \widetilde \phi_2(t)  \right] ds dt +R
\end{eqnarray}
where $\displaystyle R=R_0  - R_1 R_2 - R_1 \frac{1}{2 \pi}\int \hat g _T (t) \widetilde \phi_2 (t)dt - R_2 \frac{1}{2 \pi}\int \hat g _T (s) \widetilde \phi_1 (s)ds  $. Since $\hat g _T \in L_1(\R)$ and applying Lemma \ref{lem3} for $ j_0 =k, j_1 = k+1, j_2 =k+2, k_{gap} =l, M_1=M_2 =1$, we obtain 
\begin{eqnarray} \label{15}
	|Cov_x (a_k, a_{k+l}) | &\le &\frac{1}{(2 \pi)^2} \int \int  \left| \hat h _T (s,t) \right|    \left| \widetilde \phi (s,t) - \widetilde \phi_1(s) \widetilde \phi_2(t)  \right| ds dt + \left| R \right| \notag \\
	&\le &\frac{1}{(2 \pi)^2} \int \int  \left| \hat h _T (s,t) \right|    \left|  \phi (s,t) \varphi(s) \varphi(t) -  \phi_1(s)  \phi_2(t) \varphi(s) \varphi(t)   \right| ds dt + \left| R \right| \notag \\
	&\le & \sup_{ |s| , |t| \le \varepsilon_0}  \left|  \phi (s,t)  -  \phi_1(s)  \phi_2(t) \right| \left( \int  \left| \hat g _T (s) \right| ds \right) ^2  +|R|\notag \\
	&\le & C T^2 \kappa ^l +|R|.
\end{eqnarray}

It remains to bound of $|R|$. On the one hand, we can see that
\begin{eqnarray*}
	&\bullet& \qquad	|R_1| = \left| \E _x   [Z_k - g_T(Z_k) ] \right| =  \E _x \left| [Z_k - g_T(Z_k) ]{\bf 1}_{[|Z_k| >T]} \right| \le 2T^{-1}\E _x |Z_k|^2, \notag \\
	&\bullet& \qquad	\left| \frac{1}{2 \pi}\int \hat g _T (s) \widetilde \phi_1 (s)ds \right| =\left| \E_x g_T(Z_k) \right| \le  \E_x \left| Z_k \right| \le \E_x |a_k| + \E_x |V| \le C, \notag \\
	&\bullet& \qquad |R_2| = \left| \E _x [Z'_{k+l} - g_T(Z'_{k+l}) ] \right|  \le 2T^{-1}\E _x |Z'_{k+l}|^2, \notag \\
	&\bullet& \qquad \left| \frac{1}{2 \pi}\int \hat g _T (t) \widetilde \phi_2 (t)dt \right| = \left| \E_x g_T(Z'_{k+l}) \right|  \le \E_x |a_{l+k}| + \E_x |V'| \le C. \notag 
\end{eqnarray*}

On the other hand, similarly for $|R_0|$, we obtain
\begin{eqnarray*}
	|R_0| &=& \E_x \left[ \left|Z_k Z'_{k+l} - h_T(Z_k, Z'_{k+l} ) \right| \left( {\bf 1}_{[|Z_k| >T]} + {\bf 1}_{[|Z_k| \le T]}   \right) \left( {\bf 1}_{[|Z'_{k+l}| >T]} + {\bf 1}_{[|Z'_{k+l}| \le T]}   \right) \right] \notag \\
	&\le& \E_x \left[ \left|Z_k Z'_{k+l} - h_T(Z_k, Z'_{k+l} ) \right| \left( {\bf 1}_{[|Z_k| >T]} + {\bf 1}_{[|Z'_{k+l}| >T]}   \right) \right] \notag \\
	&\le& 2 \E_x \left| Z_k Z'_{k+l}   {\bf 1}_{[|Z_k| >T]} \right| +  2 \E_x  \left| Z_k Z'_{k+l}   {\bf 1}_{[|Z'_{k+l}| >T]} \right|  .
\end{eqnarray*}
For any positive $\delta$, let $q_\delta = {\delta+1 \over \delta}$, by Holder's inequality, we obtain
\begin{eqnarray*}
	\E_x  \left| Z_k Z'_{k+l}   {\bf 1}_{[|Z_k| >T]} \right| \le \left( \E _x |Z_k|^{2 +2 \delta}\right) ^{1 \over 2+2 \delta} \left( \E _x |Z'_{k+l}|^{2 +2 \delta}\right) ^{1 \over 2+2 \delta} \P _x(|Z_k| >T)^{1 \over q_\delta}.
\end{eqnarray*}
By Minkowski's inequality, 
\begin{align}
	\left( \E _x |Z_k|^{2 +2 \delta}\right) ^{1 \over 2+2 \delta}  & \le \left( \E _x |a_k|^{2 +2 \delta}\right) ^{1 \over 2+2 \delta} + \left( \E _x |V|^{2 +2 \delta}\right) ^{1 \over 2+2 \delta} < C, \notag \\
	\left( \E _x |Z'_{k+l}|^{2 +2 \delta}\right) ^{1 \over 2+2 \delta}  & \le \left( \E _x |a_{l+k}|^{2 +2 \delta}\right) ^{1 \over 2+2 \delta} + \left( \E _x |V'|^{2 +2 \delta}\right) ^{1 \over 2+2 \delta} < C\notag .
\end{align}
By Markov's inequality, 
\begin{align}
	\P _x (|Z_k| >T) & \le {1 \over T^{ q_\delta}} \E _x |Z_k|^{ q_\delta} \le {C \over T^{ q_\delta} }, \notag \\
	\P _x (|Z'_{k+l}| >T) & \le {1 \over T^{ q_\delta}} \E _x |Z'_{k+l}|^{ q_\delta} \le {C \over T^{ q_\delta} } \notag .
\end{align}
Hence $|R_0| \le C T^{-1}$ for $T>1$ and thus $|R| \le CT^{-1}$.

Thus, (\ref{15}) becomes $|Cov_x (a_k, a_{k+l}) | \le  C T^2 \kappa ^l + C T^{-1}$. By choosing $T= \kappa ^{- \alpha}$ with $\alpha >0$, we obtain
$$
|Cov_x (a_k, a_{k+l}) | \le C \kappa ^{l - 2 \alpha} + C \kappa ^\alpha \le C' \max \{\kappa ^{l - 2 \alpha}, \kappa ^\alpha \}.
$$
Now we choose $\alpha >0$ such that $l - 2 \alpha >0$, for example, let $\alpha = {l \over 4}$, we obtain 
$$ 
|Cov_x (a_k, a_{k+l}) | \le C \kappa ^{l \over 4}.
$$

\subsection{Proof of Lemma {\ref{lem2}}}

Inequality (\ref{e5}) follows by setting $k=l$ in (\ref{e3}) and (\ref{e4}). It suffices to prove (\ref{e4}). Recall the definition in (\ref{1.23}) and let 
\begin{eqnarray} \label{defifi}
	\psi(s) &=& \nu(P_s {\bf 1})\varphi(s) , \notag \\
	\psi(s,t; l) &=& \nu(P_s P^{l-1} P_t {\bf 1}) \varphi(s) \varphi(t) , \notag \\
	\widetilde \psi (s,t; l) &=& \psi(s,t; l) - \psi(s) \psi(t),  \\
	\widetilde \phi_0 (s,t) &=&  \widetilde \phi (s,t) - \widetilde \phi_1 (s) \widetilde \phi_2 (t) ,\notag \\
	s_{l,T}&=& \frac{1}{(2 \pi)^2} \int \int \hat h_T(s,t) \widetilde \psi (s,t; l) ds dt .\notag  
\end{eqnarray}
Then (\ref{1.28}) implies
$$
\left| Cov_x (a_k, a_{k+l}) - s_{l,T} \right| \le \left| \frac{1}{(2 \pi)^2} \int \int \hat h_T(s,t) [\widetilde \phi_0 (s,t) - \widetilde \psi (s,t; l)] ds dt \right|  + |R|.
$$
We claim that  
\begin{eqnarray}\label{cla}
	\left| \frac{1}{(2 \pi)^2} \int \int \hat h_T(s,t) [\widetilde \phi_0 (s,t) - \widetilde \psi (s,t; l)] ds dt \right|  \le C  \kappa^{k-1} T^2,
\end{eqnarray}
which implies
\begin{eqnarray} \label{12.1}
	\left| Cov_x (a_k, a_{k+l}) - s_{l,T} \right| &\le& C \kappa^{k-1} T^2 + C T^{-1},
\end{eqnarray}
which yields for any $k,m \ge 1$,
\begin{eqnarray}
	\left| Cov_x (a_k, a_{k+l}) - Cov_x (a_m, a_{m+l}) \right| \le  C  \kappa^{\min \{k-1,m-1 \} }   T^2 + C T^{-1}.
\end{eqnarray}
By choosing $T= \kappa^{-{1 \over 4}\min \{k-1,m-1 \} }$, we obtain 
\begin{eqnarray}\label{12.2}
	\left| Cov_x (a_k, a_{k+l}) - Cov_x (a_m, a_{m+l}) \right| \le  C \kappa^{\min \{{k-1 \over 4},{m-1 \over 4} \}  }.
\end{eqnarray}
Hence we can say that $( Cov_x (a_k, a_{k+l}))_l$ is a Cauchy sequence, thus it converges to some limit, denoted by $s_l(x)$. When  $k \to +\infty$,  (\ref{12.1}) becomes
$$
|s_l(x)- s_{l,T}| \le  C T^{-1}.
$$
Now let $T=T(\ell) =  \kappa^{-\ell} $, we obtain $|s_l(x)- s_{l,T(\ell)}| \le C  \kappa^\ell$. Let $\ell \to +\infty$, we can see that $s_{l,T(\ell)} \to s_l(x)$. Since $s_{l,T(\ell)}$ does not depend on $x$, so is $s_l(x)$, i.e. $s_l(x) = s_l$. Now let $m \to +\infty$ in (\ref{12.2}), we obtain
$$
|Cov_x(a_k,a_{k+l}) -s_l| \le C \kappa^{k-1 \over 4}.
$$

Now we prove the claim (\ref{cla}).
By definitions in (\ref{1.23}) and (\ref{defifi}), we obtain

\begin{eqnarray} \label{11.1}
	\left|\widetilde \phi_0 (s,t) - \widetilde \psi (s,t; l) \right| \le \left| \widetilde \phi (s,t) -  \psi (s,t; l)\right| + \left|  \widetilde \phi_1(s)  \widetilde \phi_2 (t) - \psi(s) \psi(t) \right|.
\end{eqnarray}
On the one hand, we can see that
\begin{eqnarray} \label{1.35}
	& &\left| \widetilde \phi (s,t) -  \psi (s,t; k) \right| \notag \\
	&=& \left|  P^{k-1} P_s P^{l-1} P_t {\bf 1} (x) \varphi(s) \varphi(t) - \nu (P_s P^{l-1} P_t {\bf 1} )\varphi(s) \varphi(t)\right| \notag \\
	&=&  \left|  \Pi P_s P^{l-1} P_t {\bf 1} (x) + R^{k-1} P_s \Pi P_t {\bf 1} (x) + R^{k-1} P_s R^{l-1} P_t {\bf 1} (x) -  \nu (P_s P^{l-1} P_t {\bf 1} )  \right| \notag \\
	&=& \left|  R^{k-1} P_s {\bf 1} (x) \nu(P_t {\bf 1}) +  R^{k-1} P_s R^{l-1} P_t {\bf 1} (x)  \right| \le C \kappa^{k-1}.
\end{eqnarray}
On the other hand,
\begin{eqnarray}
	\left|  \widetilde \phi_1(s)  \widetilde \phi_2 (t) - \psi(s) \psi(t) \right| &=& \left| [\widetilde \phi_1(s) - \psi(s)] \widetilde \phi_2(t) + \psi(s) [\widetilde \phi_2(t) -\psi(t)] \right| \notag \\
	&\le& \left|\widetilde \phi_1(s) - \psi(s) \right| + \left|\widetilde \phi_2(t) - \psi(t) \right| \notag \\
	&\le& \left| \phi_1(s) \varphi(s) - \psi(s) \right| + \left| \phi_2(t) \varphi(t) - \psi(t) \right| \notag ,
\end{eqnarray}
where as long as $k \ge 2$,
\begin{eqnarray*}
	\left| \phi_1(s) \varphi(s) - \psi(s) \right| &=& \left| \left[ \Pi P_s{\bf 1} (x) + R^{k-1}P_s{\bf 1} (x)  \right]\E _x [e^{isV}] - \nu(P_s {\bf 1}) \E _x [e^{isV}]\right| \notag \\
	&\le& \left| [\Pi P_s{\bf 1} (x) - \nu(P_s {\bf 1})  ] + R^{k-1}P_s{\bf 1} (x)\right| \notag \\
	&=& \left| R^{k-1}P_s{\bf 1} (x) \right| \le C \kappa ^{k-1}.
\end{eqnarray*}
Similarly, we obtain 
\begin{eqnarray} \label{1.37}
	\left|  \widetilde \phi_1(s)  \widetilde \phi_2 (t) - \psi(s) \psi(t) \right| \le C \kappa^{k-1}.
\end{eqnarray}

Therefore, (\ref{11.1}), (\ref{1.35}) and (\ref{1.37}) imply
$\left| \widetilde \phi_0(s,t) - \widetilde \psi(s,t; l) \right| \le C \kappa^{k-1} $ which yields the assertion of the claim.




\begin{thebibliography}{99}
	

\bibitem{BQ}\textsc{ Y. Benoist and J.F. Quint} (2016) Central limit theorem for linear groups,  The Annals of Probability, {\bf 44}, 1308-1340.


	
\bibitem{BL}\textsc{ Ph. Bougerol and  J. Lacroix} (1985) Products of Random Matrices with Applications to Schr\"odinger  Operators,  Birkh\"auser, Boston-Basel-Stuttgart.



\bibitem{DW} \textsc{ D. Denisov  and V. Wachtel}(2015). {\it Random walks
in cones}. arX Ann. Probab.
Volume 43, no. 3, 992--1044.

\bibitem{LPP} \textsc{ E Le Page, Marc Peigné and C Pham}(2016). {\it The survival probability of a critical multi-type branching process in i.i.d. random environment}. hal-01376091



\bibitem{FK}\textsc{ H. Furstenberg and  H. Kesten}  (1960)  Product of random matrices.
Ann. Math. Statist. {\bf 31},  457--469.  




\bibitem{GLP1} \textsc{ I. Grama, E. Lepage and M. Peigné}  (2014)  On the
rate of convergence in the weak invariance principle for dependent random
variables with applications to Markov chains. \textit{Colloquium
	Mathematicum,} {\bf 134},  1--55.


\bibitem{GLP2}\textsc{ I. Grama, E. Lepage and M. Peigné} (2016) Conditional limit theorems for products of random matrices.   Proba. Theory Relat. Field . doi:10.1007/s00440-016-0719-z


\bibitem{GH}\textsc{ Y. Guivarc'h and J. Hardy} (1988) Th\'eor\`emes limites pour une classe de cha\^ines de Markov et applications aux diff\'eomorphismes d'Anosov. Ann. Inst. H. Poincar\'e  Proba. Statist.  {\bf 24},     Proba no. {\bf 1}, 73--98.


\bibitem{GR}\textsc{ Y. Guivarc'h and  A. Raugi} (1985) Frontiere de Furstenberg, propri\'et\'es de contraction et theorems de convergence. Z. Wahrsch. Verw. Gebiete, {\bf 69}, 187-242.





\bibitem{H2}\textsc{H. Hennion} (1997) Limit theorems for products of positive random matrices, Ann.  Proba. {\bf 25}, no. {\bf 4},  1545--1587. 

 \bibitem{H1}\textsc{H. Hennion} (1984) Loi des grands nombres et pertubations pour des produits reductibles de matrices al\'eatoires ind\'ependantes, Z. Wahrsch. Verw. Gebiete, {\bf 67}, 265-278.
 
 \bibitem{HH}\textsc{H. Hennion and  L. Hervé} (1984) Stable laws and products of positive random matrices,Journal of Theoretical Probability {\bf 21}, no. {\bf 4},  966--981. 
 
\bibitem{LePage82} { \sc E. Le\ Page}  (1982)  Th\'{e}or\`{e}mes limites pour les
produits de matrices al\'{e}atoires, Springer Lecture Notes
{\bf 928}, 258--303.

\bibitem{Levy} { \sc P. Levy }  (1937)  Th\'{e}ory de l'addition des variables aléatoires. Gauthier-Villars 


\end{thebibliography}
\end{document}